\documentclass{article}
\usepackage{latexsym}
\usepackage{epsfig}
\begin{document}

\title{\huge A Note About Universality Theorem as an Enumerative 
 Riemann-Roch Theorem}
\author{Ai-Ko Liu\footnote{email address: 
akliu@math.berkeley.edu} \footnote{Current Address: 
Mathematics Department of U.C. Berkeley}\footnote{
 HomePage:math.berkeley.edu/$\sim$akliu}}
%\date{April, 8}
\maketitle
%\centerline{\bf Research Announcement}
%\centerline{\bf Preliminary Version }
\newtheorem{main}{Main Theorem}
\newtheorem{theo}{Theorem}
\newtheorem{lemm}{Lemma}
\newtheorem{prop}{Proposition}
\newtheorem{rem}{Remark}
\newtheorem{cor}{Corollary}
\newtheorem{mem}{Examples}
\newtheorem{defin}{Definition}
\newtheorem{axiom}{Axiom}
\newtheorem{conj}{Conjecture}

This short note is a supplement of the longer paper [Liu6], in which the
 author gives an
algebraic proof of the following universality theorem.

\begin{theo}\label{theo;  main1}
Let $\delta\in {\bf N}$ denote 
 the number of nodal singularities.
Let $L$ be a $5\delta-1$ very-ample\footnote{For the
definition of $k$-very ample line bundles, please consult [Got]. The
$5\delta-1$ power of very ample line bundles are $5\delta-1$ very ample.}
 line bundle on an algebraic surface $M$,
 then the number of $\delta-$nodes nodal singular curves in a generic 
 $\delta$ dimensional linear sub-system of $|L|$ can be expressed as 
 a universal polynomial (independent to $M$) 
of $c_1(L)^2$, $c_1(L)\cdot c_1(M)$, $c_1(M)^2$, 
$c_2(M)$ of degree $\delta$.
\end{theo}

\medskip

 The finiteness of the 
``number of $\delta-$nodes nodal singular curves in a generic 
 $\delta$ dimensional linear sub-system of $|L|$'' was proved by
 G${\ddot o}$ttsche in [Got] proposition 5.2. Our theorem shows that these
 numbers are topological invariants of $(M, L)$.

 A weaker form of the above statement has appeared in [Got] based on
 results in [V] and [KP] for small $\delta$. In his conjecture,
 G${\ddot o}$ttsche had assumed the existence of a lower bound $m_0$, such that
 for any very ample $L_0$, the statement in the above theorem 
 holds for $L=L_0^{\otimes m}$, with $m\geq m_0$.

The purpose of the
note is to address upon the geometric and 
topological meanings of the universality theorem
and its relationship with
 the well known surface Riemann-Roch formula of algebraic surfaces,
 to compare the difference of 
the symplectic and the algebraic approaches in resolving this problem, and 
 list the open problems and conjectures 
related to the solution of the problem, etc. in a less technical term.

\bigskip

  The universality theorem has played an important role in giving an in-depth
 understanding of Yau-Zaslow conjecture [YZ] (see page \pageref{YauZaslow}
 also).  
  The universality theorem and the various machineries built up in
  [Liu1], [Liu3], [Liu4], [Liu5], [Liu6] and [Liu7] are used in [Liu2]
 to solve the Harvey-Moore conjecture [HM] on the counting rational curves of
 Calabi-Yau K3 fibrations. The generalization of our universality theorem
 to higher dimensions, counting divisors in a complete linear system, also
 bring us new machineries to understand the concept of ``number of 
nodal curves'' in higher dimensions.

\medskip

 In section \ref{section; meaning} we outline the formulation of the 
general enumerative problem of a complete linear system 
which relates the universality theorem with the 
classical Riemann-Roch theorem. In fact the above
 universality theorem can be interpreted
naturally as the prolongation of the surface Riemann-Roch formula.

In section \ref{section; outline}, we outline the key ideas of the proof in
 [Liu6] and address the key conceptual and technical 
issues the paper has resolved. As the
 algebraic proof in [Liu6] is slightly lengthy, we hope that the
 sketch of the basic ideas may help the reader to digest the long paper.

In section \ref{section; open}, we list a few open problems related to the
 solution of universality theorem.

\medskip

\section{The Geometric Meaning of The Universality Theorem}\label{section; meaning}

\medskip

 Recall the well known classical surface Riemann-Roch formula (e.g. [GH] or
 [Z]) that given  an algebraic surface $M$, the
 Euler characteristic of an algebraic line
 bundle $L$ on $M$ satisfies, \label{RR}

$$\chi(L)\equiv h^0(M, L)-h^1(M, L)+h^2(M, L)=\chi({\cal O}_M)+
 {c_1(L)^2-c_1(K_M)\cdot c_1(L)\over 2}$$
$$={c_1(M)^2+c_2(M)\over 12}+{c_1(L)^2-c_1(K_M)\cdot c_1(L)\over 2}.$$

When $L={\cal O}_M$, it is reduced to the well known Noether formula
 $\chi({\cal O}_M)={c_1(M)^2+c_2(M)\over 12}$.

 As a special case of the general Hirzebruch-Riemann-Roch formula [Hir],
 the above formula can be viewed as the grade two term of the expansion of 
 $Todd(M)ch(L)$, where $Todd(M)$ and $ch(L)$ are the well known Todd class 
 and Chern Character (please consult [F] or [Hir] for their definitions).

 On the other hand, a more 
down-to-earth approach to the surface Riemann-Roch formula 
involves showing that the holomorphic Euler number
$\chi(L)$ can be expressed as a universal degree one polynomial of
 $c_1(M)^2$, $c_2(M)$, $c_1(M)\cdot c_1(L)$ and $c_1(L)^2$.  I.e.
$$\chi(L)=A_1c_1(M)^2+A_2c_2(M)+A_3c_1(M)\cdot c_1(L)+A_4c_1(L)^2.$$

The above statement can be viewed as an existence result of the universal
 polynomial in $c_1(M)^2$, $c_2(M)$, $c_1(L)^2$, $c_1(L)\cdot c_1(M)$ to express
 the ``topological'' quantity $\chi(L)$.
 Once such a polynomial is shown to exist,
then we identify the four coefficients $A_1, A_2, A_3, A_4$ 
through concrete examples.
While there are plenty of choices of algebraic surfaces 
in identifying these coefficients 
$A_1, A_2, A_3$, and $A_4$, a natural choice to identify $A_1$ and $A_2$ (i.e. 
 the coefficients in the Noether formula)
 is to consider $M=K3$ and $M={\bf P}^2$. Then by a simple arithmetic and 
 the knowledge about their irregularities $q$ ($={b_1\over 2}$, geometric
 genus $p_g$ 
($={b^+_2-1\over 2}$), we find 
 $\chi({\cal O}_{K3})=1-0+1=2=A_1\times 0+A_2\times 24$,  so $A_2={1\over 12}$. 
Similarly
 $\chi({\cal O}_{{\bf P}^2})=1-0+0=A_1 \times c_1({\bf P}^2)^2+A_2c_2({\bf P})=
9A_1+3A_2$, so $A_1={1\over 12}$, too.  Then one may identify 
 $A_3$ and $A_4$ by adopting the polarized manifolds $({\bf P}^2, {\bf H}^d)$ 
 and $(T^4, L)$. By the classical calculation
$h^0({\bf P}^2, {\bf H}^d)={d(d+3)\over 2}$, 
 $h^1({\bf P}^2, {\bf H}^d)=h^2({\bf P}^2, {\bf H}^d)=0$ and
 $h^0(T^4, L)={c_1(L)^2\over 2}$, 
$h^1(T^4, L)=h^2(T^4, L)=0$ for any ample $L$, the coefficients 
$A_3$ and $A_4$ can be identified with ${1\over 2}$.

In this way we have recovered the Riemann-Roch formula,

$$\chi(L)={1\over 12}(c_1(M)^2+c_2(M))+{1\over 2}(c_1(L)^2-c_1(L)\cdot c_1(K_M)).$$

\begin{rem}\label{rem; comment}
The choices of ${\bf P}^2, ({\bf P}^2, {\bf H}^d), K3$ and $(T^4, L)$
 can be justified by the following facts, 

(i). The dimension formulae of $h^i({\bf P}^2, {\bf H}^d)$
 are well known classically (Consult for example [Ha] chapter III.5). 

(ii). The $K3$ and $T^4$ have trivial canonical
 bundles \footnote{therefore $h^i(L)=h^{2-i}(L^{\ast})$ under 
 Serre duality.} and therefore have vanishing first Chern classes.
\end{rem}

\medskip

\begin{rem}\label{rem; identify}
The identification of the second term 
${c_1(L)^2-c_1(L)\cdot c_1(K_M)\over 2}$ with the
expected complex dimension, Gromov-Taubes formula ${C^2-C\cdot c_1(K_M)\over 2}$, 
(with the substitution $c_1(L)=C$) of moduli space of curves
 is the starting point of our projects to formulate
the algebraic Seiberg-Witten theory exploring the
correspondence of surface Riemann-Roch theorem
and Taubes' version of Gromov invariant. (see [Liu1], [Liu3]
 and [Liu6]).
\end{rem}

\bigskip

\subsection{Linear Systems and Singular Curves}\label{subsection; formulate}

\bigskip

 Given a line bundle $L\mapsto M$ over the algebraic surface $M$,
 the complete linear system
 ${|L|=\bf P}(H^0(M, L))$ parametrizes the effective curves determined by
 the non-zero global sections $\in H^0(M, L)$. So we may identify 
 $|L|$ with a projective space parametrizing the linear equivalence classes of 
effective divisors.
When $L$ is not ``positive''
 enough, the dimension $h^0(L)$ depends on the algebraic surface $M$ and $L$ 
explicitly and has to be investigated individually. On the other hand for
 sufficiently positive \footnote{for example, when 
 $L\otimes K_M^{-1}$ is ample.} $L$ the Kodaira vanishing theorem implies 
$h^1(L)=h^2(L)=0$. 
 Then $h^0(L)=h^0(L)-h^1(L)+h^2(L)=\chi(L)$ is a topological quantity (i.e.
 independent to the complex/holomorphic structures 
of $(M, L)$) and is determined by
 surface Riemann-Roch formula. 

 Let $L$ satisfies the condition that $L^{\ast}\otimes K_M$ is non-nef.
 Under such a simplifying condition \footnote{This condition appears naturally
 in the definition of algebraic family Seiberg-Witten invariant.} 
the $h^2(L)=h^0(L^{\ast}\otimes K_M)=0$ by Serre duality and the fact
 $L^{\ast}\otimes K_M$ has a negative degree.
 The primary invariant one may associate with the projective space 
$|L|$ is its expected dimension $h^0(L)-h^1(L)-1=\chi(L)-1$, 
computed by the Riemann-Roch formula.

 On the other hand, the well known classical Bertini theorem (see e.g.
 page 179 of [Ha]) implies that for any very ample (or based point free) 
$L$, the generic
 divisors \footnote{They are curves when $dim_{\bf C}M=2$. In the paper we
 use the term ``curves'' and ``divisors'' interchangeably.} 
in $|L|$ are smooth and irreducible. Moreover the curves
 tend to develop singularities when they move to the boundary points of the
 top stratum of smooth curves. 
This suggests us to stratify the space
 $|L|$ into the various strata according to the ``topological types'' 
\footnote{The topological type here refers to the topological type of the
 link space in ${\bf S}^3\subset {\bf C}^2$ of the isolated singularity.} of the 
list of singularities of the divisors (curves) and the unique 
top dimensional open stratum
parametrizes smooth curves in $|L|$.

\bigskip 

The intersection pairing $\int_{|L|}{\bf H}^{\chi(L)-1}=1$ can be
 interpreted as, 

\medskip

\noindent (i). The intersection number associated with the following enumerative
 problem: How many irreducible smooth curves in $|L|$ are there which
 pass through $\chi(L)-1$ generic points in $M$?  

or equivalently,

\noindent (ii). 
The number of smooth irreducible curves in a generic $0$ dimensional
 linear sub-system of $|L|$. 

\medskip

 Even though the sub-scheme of $|L|$ parametrizing 
smooth curves through $\chi(L)-1$ generic
 points in $M$ may not be reduced and smooth of zero dimension, 
the intersection
 number defined by intersection theory [F] is still defined and
 is equal to $1$.

 If one is familiar with the definition of algebraic Seiberg-Witten invariant
 of $L$, it will be clear that 
 the above number is nothing but the algebraic Seiberg-Witten invariant
\footnote{The $t_L\in T(M)$ insertion in the $SW$ invariant 
indicates that we have restrict the
holomorphic structure of ${\cal E}_{c_1(L)}$ to $L$.}  of
 $L$, ${\cal ASW}(t_L , c_1(L))$.

\medskip

\noindent {\bf The Relationship with A Degenerated
 Case of Donaldson-Thomas Invariant}

\medskip

 Consider a closed sub-scheme $Z$ of a smooth algebraic scheme $M$, 
 there is a associated short exact sequence of the ideal sheaf 
 ${\cal I}_Z$,

 $$0\mapsto {\cal I}_Z\mapsto {\cal O}_M\mapsto {\cal O}_Z\mapsto 0.$$

The Hilbert scheme is the universal object parametrizing the closed sub-schemes
 (or equivalently ideal sheaves) with a fixed Hilbert polynomial.

 In [DT], [Th], S. Donaldson and R. Thomas considered the moduli space
 of stable sheaves on Fano or Calabi-Yau $3-$fold $M$ 
and defined a version of ${\bf Z}$ valued invariant
based on the virtual fundamental class technique of [LiT1]. When one
 specifies the sheaf rank to $1$, the general problem collapses to the
 definition of the virtual fundamental cycle of Hilbert schemes. 

 In the same paper [Th], R. Thomas also considered the moduli space
 of ideal sheaves ${\cal I}_Z$ with a fixed determinant \footnote{The determinant
 of a coherent sheaf is defined by taking the alternating product of
 determinants of a locally free resolution of ${\cal I}_Z$.} and defined
 the Donaldson-Thomas invariant similarly. 

 These Donaldson-Thomas invariant
 was compared in a recent paper [MNOP] with the Gromov-Witten
invariants for the ``local Calabi-Yau'' in relating the Donaldson-Thomas invariants
 with the mysterious ``Gupakurma-Vafa invariants''.

  When $dim_{\bf C}M=2$ and we consider the
 codimension one $Z\subset M$, the ideal sheaf ${\cal I}_Z$ becomes
 invertible. In such a case the determinant of ${\cal I}_Z$ is 
${\cal I}_Z\cong {\cal O}_M(-Z)$ 
itself and fixing the determinant corresponds to fixing the linear
 equivalence class of $Z$.

 Therefore the complete linear system $|L|$ is nothing but the
 $dim_{\bf C}M=2$ 
analogue of Donaldson-Thomas moduli space of ideal sheaves with a fixed
determinant. $\Box$
\label{DT} 

\subsection{\bf Enumeration of Singular Curves} \label{subsection; singular}

\bigskip

 On the other hand, there is no reason to restrict ourself to
 enumerating smooth curves in $|L|$ only.
 Instead one may consider all the (non-compact) strata of singular curves
 parametrized by the list of topological types of 
singularities. 

 It is known that the creation of isolated singularities on an algebraic
 curve drops
 the ``expected dimension'' of the moduli space of 
curves and the geometric genera of the curves. E.g. in the case of nodal 
singularities, introducing a single node on the curve 
drops the ``expected dimension'' and the genus
 of the curve by $1$.

In this section we let
 the bold character ${\bf m}$ denote a finite list of isolated curve
 singularities in algebraic surfaces \footnote{Because our focus is to give
a conceptual understanding, we will go into the details of how to
 express ${\bf m}$ in terms of the topology and the singular multiplicities
 of the list of singularities.}.

\medskip

\noindent {\bf Question 1}: Let $d_{\bf m}(L)$ denote the ``expected dimension'' of
 the stratum (which can contain more than one component) 
of singular curves with the prescribed ``topological types'' of
 singularities. How many singular curves are there which pass through the
 generic $d_{\bf m}(L)-$number of points? Or how many singular curves are there 
 in the closure of the 
stratum which lie in the generic $\chi(L)-1-d_{\bf m}(L)$ dimensional
 linear sub-system of $|L|$?

 At first the question may look ill-posted as 

\noindent (a). The stratum of singular curves (fixing the topological type)
  is usually non-compact. 

It can be remedied by adding the 
compactifying strata to the given stratum \footnote{this means
 adding strata corresponding to singularities into which the original
 singularities can degenerate}. One may compactify it inside $|L|$. 

\medskip

\noindent (b). The sub-scheme of the union of compactified strata 
parametrizing the ${\bf m}$ singular 
curves in the generic $dim_{\bf C}|L|-d_{\bf m}(L)=
\chi(L)-1-d_{\bf m}(L)$ dimensional
 linear sub-system is ``expected'' to be regular of zero dimension if
 all the generic conditions are met. But it may not be regular of zero dimension
 geometrically.

 Nevertheless the well known remedy is to
 rephrase the above question in terms of 
intersection theory technique and ask:

\medskip

\noindent {\bf Question 2}: How can we attach an intersection number (usually
 called the virtual number) to each compactified 
stratum of singular curves in $|L|$ with prescribed ``topological type
 of singularities''\footnote{We fix the topological types
instead of the analytic types of the germ of singularities here. If we work
 over ${\bf C}$, the topological type of isolated curve singularities is
 determined by the link space of the singularity in ${\bf S}^3\subset {\bf C}^2$,
 which is a multi-component link in ${\bf S}^3$.}
such that it is reduced to the actual curve counting when the sub-scheme 
of curves in the generic $\chi(L)-1-d_{\bf m}(L)$ 
dimensional linear sub-system of $|L|$ is  
regular of the expected dimension $0$?

In the following, we will abbreviate and denote these numbers $n_{\bf m}$ 
as ``the number of singular curves''
 in $|L|$.  It should be understood in the sense of virtual numbers unless
 we are capable of proving a version of regularity result on the sub-scheme 
of of $|L|$ parametrizing the singular curves
 in generic $dim_{\bf C}|L|-d_{\bf m}(L)$ dimensional sub-linear system.

\medskip 

  It is highly non-trivial to check whether these 
``virtual numbers of singular
 curves'' are invariants. It makes perfect sense to ask:

\medskip

\noindent {\bf Question 3}: Are the ``number of singular curves'' well
 defined?
Suppose that the ``number of singular curves'' 
has been defined, how does the intersection number
 depend on the algebraic surfaces $M$ and the line bundle $L$? Are they
 deformation invariants?
 
 If there is a family of algebraic surfaces with sufficiently 
very ample relative polarizations which restrict to $L$ on the special fiber, are
 the corresponding ``numbers of singular curves''
 defined for each member within such a family independent
 to the deformation?  

\bigskip

 Are these numbers 
 ``topological'', which depend only on the topological data on $(M, L)$?
As the dimension $dim_{\bf C}|L|=\chi(L)-1$ is topological
 (by surface Riemann-Roch) when $L$ is sufficiently positive, would it be too
 much to expect those ``virtual numbers of singular curves'' do, too?

\medskip

 Before we investigate the relationship of 
our universality theorem with these structural
questions on the ``virtual number of singular curves'', 
let us review how the special cases of these questions have
shown up in the researches of the various group of people and how they have
answered them.

\bigskip

Among the different type of isolated curve singularities, 
the nodal singularities, i.e. the normal crossing curve singularity, is the
 most and well-studied one.

Take $M={\bf CP}^2$ and $L={\bf H}^d$, 
the varieties of irreducible nodal curves in $|L|=|{\bf H}^d|$ over 
${\bf CP}^2$ are known
 to be Severi varieties. In [H] it was proved that the Severi varieties
 are irreducible. On the other hand,
the number of
 $\delta$-node nodal curves in the generic $\delta$ dimensional linear
 subsystem of $|{\bf H}^d|$ is also known to be
 the Severi degree of the Severi variety and was 
calculated by [CP] (see also [Ran]) in terms of recursive formulae on
 the degree $d=c_1({\bf H}^d)$. 
The same authors had shown that
 the intersection numbers could be understood in the
 classical sense, i.e. the ``number of nodal curves'' really counts the 
 discrete number of nodal singular curves in $|{\bf H}^d|$. 

On the other hand, when \footnote{I.e. $c_1(M)=0$ and $\pi_1(M)=\{1\}$.}
 $M=K3$ and $c_1(L)\in H^2(M, {\bf Z})$, the
 ``numbers of nodal curves'' usually cannot be understood 
 in the classical sense. This is the perfect examples to observe that
 the stratum of $|L|$ parametrizing the $\delta-$ node nodal curves can be
 rather ill-behaved for specific choice of algebraic K3 and $L$.

 In particular, when $\delta={c_1(L)^2\over 2}-1$,
 these ``number of rational nodal curves''
 attached to all such $|L|$ are predicted by
 the Yau-Zaslow formula [YZ] to be generated by the modular form 
 ${1\over \eta(q)^{24}}$ (the generating
 function of the Euler numbers of the Hilbert Schemes of $K3$, by a beautiful
 result
 of G${\ddot o}$ttsche).

  Recall the conjecture of Yau-Zaslow makes the following three 
predictions \footnote{The (3). below is usually known to be the
 Yau-Zaslow conjecture sometimes. To understand why the (1). and
 (2). are an integral part of the conjecture, please consult the
 original argument in the paper of Yau and Zaslow [YZ].}.

\medskip

\noindent (1). The ``number of embedded nodal rational curves'' in $|L|$ 
is well defined. \label{YauZaslow}

\medskip

\noindent (2). The ``number of embedded nodal rational curves'' in $|L|$
 depends on $L$ only \footnote{This is a very strong statement which
 claims that the ``number of embedded nodal rational curves'' does not
 depend on the choices of the classes $c_1(L)\in H^2(M, {\bf Z})$ if they
 have the same self-intersection number.}
through $\delta={c_1(L)^2\over 2}-1$, or equivalently, the self-intersection
number of $c_1(L)$, and can be 
denoted by $n_{\delta}$.

\medskip

\noindent (3). The generating function $\sum_{\delta\geq 0}n_{\delta}q^{\delta}$
 coincides with $\bigl({1\over \prod_{i=0}^{i=\infty}(1-q^i)}\bigr)^{\chi(M)}$.

Due to the difficulty in interpreting
 the ``number of nodal rational curves'' classically, some people (e.g. [BL])
 had attempted
 to interpret these numbers as the Gromov-Witten invariants of the ${\bf S}^2$
 families in some
 special cases. Indeed they
were able to justify the (iii) statement in the 
conjecture (replacing the
 term ``number of nodal rational curves'' by genus zero Gromov-Witten
 invariants) of ${\bf S^2}$ hyperkahler families 
when $L$ is primitive, i.e. $c_1(L)$ 
 being a primitive element in the $K3$ lattice $H^2(K3, {\bf Z})$. 

  Despite that Gromov-Witten invariants on algebraic surfaces 
are closed related to the concept of ``number of nodal curves'' we 
proposed above and it may be tempting
 to think Gromov-Witten invariants as the exact intersection 
theoretical interpretation of {\bf Question 2}, we point out the 
importance to separate these two different concepts. Especially it becomes
 apparent in some special case on K3 [Gat]
that the intersection numbers predicted by 
Yau-Zaslow formula differs from the genus zero ${\bf S}^2$ family 
Gromov-Witten invariant of the
 same (non-primitive) class. On the other hand, Yau-Zaslows' prediction matches up 
perfectly with calculation by Vainsencher [V] and was the strong
 evidence in the early days that Yau-Zaslow formula gave the correct
prediction (See the final section ``conclusions and Prospects'' of [YZ]).

A careful
investigation upon their subtle difference shows us that genus zero 
Gromov-Witten invariant
enumerates not only maps onto immersed curves but also 
the multiple-coverings of maps. Yet the ``number of nodal curves''
in Yau-Zaslow conjecture counts immersed nodal curves in the linear
 system.
 Only when the class is primitive 
\footnote{and the curve cone/Picard lattice is one dimensional.}, 
these two concepts become coincide accidentally--as there is no chance of having 
multiple-covering maps into any primitive class.
 
\medskip

 On the other hand, besides some special cases of cohomology classes
 and almost complex structures closed enough to integrable complex
 structures, the author is not aware of the
 general definition of ``number of nodal curves'' in symplectic geometry.

\medskip

\begin{rem}\label{rem; relate}
The above discussion suggests that we should consider the ``number of
 nodal curves'' and Gromov-Witten invariants as separated concepts. They
 coincide only in special cases.

 In fact the definition of Gromov-Witten invariants depend on the moduli stack
 of marked curves, and is related to the ``world-sheet'' point of view of
 string theory.  Our theory of universality theorem makes use of the
 universal spaces of algebraic surfaces, which is closely related to the 
 space-time approach of string theory.
 Nevertheless, we expect that there should be a
 combinatorial formula relating these two types of concepts.
 Our algebraic construction of the algebraic family obstruction bundle
 \footnote{embedded in the proof of the main theorem.} (along with the technique
 to remove the type $II$ exceptional class 
contributions when $L$ is not sufficiently positive) 
provides an algebraic way to
define the ``number of nodal curves'' as a virtual number, even when the
 corresponding moduli space of nodal curves is not well behaved.
\end{rem}

 For algebraic surfaces with geometric genera $p_g>0$, 
Gromov-Witten invariants are frequently vanishing (due to the simple type 
condition. See [Liu3] section 4.3), while 
``the number of nodal curves'' in the linear system $|L|$ can still be
 very rich! In the case of $K3$ (with $p_g=1$) 
it is also the reason why people had adopted the twistor families
 of $K3$ in counting curves.

 Some examples in the table I of [Va] indicates
 that the ``number of nodal curves''
 may not always be topological--independent to degeneration of complex
structures.

\medskip

\noindent{\bf Topological or Non-topological?}

\bigskip

Does the above discussion
 mean that we have to give up the idea that the ``numbers
 of nodal curves'' are topological?  Not really!

\medskip

 In fact the real essence of the universality theorem is to tell us that
 for the number of nodal singularities $=\delta$:

\medskip

\noindent
 (1). When $L$ is $5\delta-1$-very ample, the ``number of nodal curves'' on
 any algebraic surface  
 can be understood in the classical sense, once we take into account of
 multiplicities. The argument is essentially due to G${\ddot o}$ttsche [Got].

\medskip

 \noindent (2). These numbers are actually topological!  Indeed they can be
expressed as degree $\delta$ 
universal polynomials of $c_1^2(L), c_1(L)\cdot c_1(M)$, 
$c_1^2(M)$, $c_2(M)$ in the same way that the expected dimension 
$dim_{\bf C}|L|$ has been given
 by a degree one polynomial of $c_1^2(L), c_1(L)\cdot c_1(M)$, 
$c_1^2(M)$, $c_2(M)$ through the Riemann-Roch formula (see page \pageref{RR}). 

\medskip

 \noindent (3). The theorem implies that the ``number of nodal curves'' on totally
 different algebraic surfaces are intimated related, even though the
 algebraic geometry on these distinct surfaces can be quite different. The theorem
 allows us to analyze the dependence of these enumerative geometric 
information upon the underlying surface $M$ and $L$
 in terms of homotopic types of the algebraic 
surfaces $M$ and the cohomology class $c_1(L)$ of $L$.

\medskip

On the other hand, we address the following apparent puzzles briefly.

\medskip

\noindent {\bf Question}: If the universality theorem asserts that
 the counting of nodal curves is topological, why are there examples
 of deformation equivalent surfaces (check Table 1 on page 11 of 
[Va]) which have different
 ``numbers of nodal curves'' in the same classes?  

 In a simple term to answer, it is due to
 the $5\delta-1$-very ample condition
 imposed on $L$.  We may rephrase the above question by the following one,

\medskip

\noindent {\bf Question}: Does the universality theorem mean that
the nodal curve counting (or more generally singular curve 
counting) is completely rigid
 as the ``number of singular curves'' are always predicted by one single 
universal formula?

\medskip

 Not really!  Each algebraic surface still preserves its own character. 
 The universality theorem only implies that in a complete linear system $|L|$
 with $L$ far away from the origin of
the ample cone, the ``number of nodal curves'' behaves
 topologically. Our universality theorem gives the effective bound
 on $L$ for such statements to hold.

On the other hand, without the suitable very ample condition, 
the topological prediction extrapolated
from the universality theorem is generally not the exact answer 
\footnote{This comment does not apply when $M=K3$ or $T^4$ in which case
the topological answer by the universality theorem matches with the
geometric answer perfectly despite that $L$ may not be $5\delta-1$ very ample.
This is because a simple vanishing argument on the contribution from 
 type $II$ exceptional curves. See page \pageref{subsubsection; vanishing}.}!

 Then we may ask:

\noindent {\bf Question}: What happens if $L$ is not $5\delta-1$ very ample?

 For \underline{general}
 algebraic surfaces, apparently the predictions on the
 ``numbers of nodal curves'' are not accurate without $5\delta-1$ very ample
condition. 

 For the Calabi-Yau algebraic surfaces, i.e. K3 surfaces or Abelian surfaces
 ($\cong T^4$), the prediction from the universal polynomials of our
 universality theorem still holds even without the positivity condition on $L$.
  The reason that $K3$ and $T^4$ are special in this aspect has to be 
understood in a rather theoretical level. We will discuss the 
vanishing result in subsection \ref{subsubsection; vanishing} briefly, following
 the spirit of Riemann-Roch theorem.

 On the other hand, the fact that for $\delta<7$ the Zaslow-Yau prediction
and Vainsencher's enumeration of the universal polynomials {\bf Do} match has 
been a convincing evidence in the early days when Yau-Zaslow had made their
 conjecture [YZ].

\subsection{\bf The Type II Exceptional Curves and The
 Discrepancy to Universality Theorem}

\bigskip
 
To understand the non-topological nature of the ``numbers of
 nodal curves'' generally 
(without any assumption on $L$) and how the geometric answers derivate from
 the universality theorem,
 we have to reflect the way we deal with  
the similar phenomenon in
 the original Riemann-Roch theorem. While it is true that for
 $L\otimes K_M^{-1}$ ample, $h^1(L)=h^2(L)=0$ the $h^0(L)=\chi(L)$ and
 therefore $dim_{\bf C}|L|$ is topological, generally the expected dimension
 $h^0(L)-h^1(L)-1$ is not topological and
 can not be predicted by a topological formula.

 If we go back to the history of the development of Riemann-Roch
 theorem (see e.g. chapter VII-2A of [Di]/or chapter IV and its appendix of 
[Za]), the higher sheaf
cohomologies were introduced (in the classical language, they were known to
be $h^1(M, L)=$ the super-abundance and $h^2(M, L)=$ index of specialty)
 exactly to balance the discrepancy between Riemann-Roch
 formula (of topological nature involving characteristic 
classes on $M$ and $L$) and $h^0(L)$.

  If one browses through the paper [Liu3] for the construction of
 the algebraic (family) Seiberg-Witten invariants, it is clear that
 the construction of the algebraic Seiberg-Witten invariant has been
 separated into three cases.

\medskip

\noindent (i). The case of regular linear system with $h^1(L)=h^2(L)=0$: 
In this case, the space $|L|$ is smooth of 
the right dimension $h^0(L)-1=\chi(L)-1$, and the algebraic Seiberg-Witten
 invariant can be defined in the classical sense.

\medskip

\noindent (ii). $h^2(L)=0$ but $h^1(L)$ may be non-zero: In this case
 $|L|$ is not of its expected dimension. Yet one may introduce
 the algebraic Kuranishi model\footnote{Read [Liu3] for more details.}
 of bundle maps 
$\Phi_{{\bf V}{\bf W}}:{\bf V}\mapsto {\bf W}$ such that
 $rank_{\bf C}{\bf V}-rank_{\bf C}{\bf W}=h^0(L)-h^1(L)=\chi(L)$. 

 Accordingly the algebraic Seiberg-Witten invariant is defined to be
 $\int_{{\bf P}({\bf V})} c_{top}({\bf H}\otimes {\bf W})$.

\medskip

\noindent (iii). $h^2(L)\not=0$: In this case $h^2(L)$ is the
 correction term to
 the dimension of the class $c_1(L)$ from the Riemann-Roch formula.

 This suggests us to group $h^0(L)-h^1(L)$ inside $h^0(L)-h^1(L)+h^2(L)$ and
the expected dimension of the linear system $|L|$ in this case is 
 $h^0(L)-h^1(L)-1=\chi(L)-1-h^2(L)$. Namely, we have to calculate
 $h^2(L)$ and subtract it from $\chi(L)-1$ to get the actual ``expected
 dimension'' $h^0(L)-h^1(L)-1$ and $h^2(L)$ plays the role of the
 ``correction term'' in the following formula,

$$h^0(L)-h^1(L)-1=\{\chi(L)-1\}-h^2(L).$$

Then we may ask

{\bf Question}: How does the correction term $h^2(L)$ of the surface
Riemann-Roch formula shed light on the discrepancy of ``number of nodal curves''
 to the Universality theorem? 

Is there any analogous concept
 in the enumeration theory of nodal curves parallel to the dimension of
 second sheaf cohomology 
 $h^2(L)$ in the calculation of the expected dimension $h^0(L)-h^1(L)-1$?

 In fact, such an analogue exists.
It is the excess contribution to the family invariants 
from decomposition of curves involving type $II$ exceptional curves. It plays
an exact analogue of $h^2(L)$ in the surface Riemann-Roch formula.

\begin{defin} \label{defin; type2}
 A divisor class $e$ is an exceptional class over an algebraic fibration
 ${\cal X}\mapsto B$ if it satisfies the following condition.

 (i). The self intersection number $e^2<0$.

 (ii). The degree of $e$ with respect to ample polarizations is positive.
\end{defin}

 The first condition ensures that any irreducible representatives of $e$
  is unique in the corresponding fiber. The second condition is necessary for
 $e$ to be represented by algebraic curves. The class $e$ is exceptional
 in the sense that for an irreducible ${\bf e}$ representing $e$, 
$h^0({\bf e}, {\cal O}({\bf e}))=0$ and the curve ${\bf e}$ is infinitesimally
 rigid.

 In the family theory approach to the universality theorem, we encounter
 two types of exceptional curves on the universal families $M_{\delta+1}\mapsto 
M_{\delta}$ 
and they are called type $I$ and type $II$ exceptional curves, respectively.

 The fibration $M_{\delta+1}\mapsto M_{\delta}$ can be constructed from the product
 family $M\times M_{\delta}\mapsto M_{\delta}$ by blowing up $n$-consecutive cross sections
 of the intermediate blown up fiber bundles\footnote{Read the beginning of 
section 2 [Liu6] for
 more details.}. So there is a canonical projection map $M_{\delta+1}\mapsto M\times 
 M_{\delta}\mapsto M$. Following the convention in [Liu6],
 the exceptional divisors of the intermediate blown up manifolds are
 denoted by $E_1, E_2, \cdots, E_{\delta}$, respectively. 

\begin{defin}\label{defin; one}
A type $I$ exceptional classes is an exceptional class
 of the form $e=E_i-\sum_{j_i>i}E_{j_i}$.
\end{defin}

 It lies
in the kernel of the proper push-forward 
${\cal A}_{2\delta+1}(M_{\delta+1})\mapsto {\cal A}_1(M)$.

 An irreducible algebraic curve representing the type $I$ exceptional class 
 is said to be a type $I$ exceptional curve in the family $M_{\delta+1}\mapsto M_{\delta}$.

 The type $I$ exceptional curves have played important roles in the proof 
of universality theorem.

\begin{defin}\label{defin; two}
A type $II$ exceptional class in the universal family $M_{\delta+1}\mapsto M_{\delta}$ 
is an exceptional class which is not in the kernel of the proper push-forward
 ${\cal A}_{2\delta+1}(M_{\delta+1})\mapsto {\cal A}_1(M)$.
\end{defin}

 Over the universal families $M_{\delta+1}\mapsto M_{\delta}$ 
the appearance of the type $II$ exceptional
 classes and their excess contributions to the family invariants 
 can be understood from three different prospectives,

\medskip

\noindent (a). From the linear system interpretation of the nodal curve counting,
there are strata in $|L|$ which consist of curves with non-isolated
singularities. These are the curves with non-reduced irreducible
components. Just like the counting of nodal curves using the universal
 family may encounter other type of singular curves with isolated
 singularities, sometimes
 the singular curves with non-reduced irreducible components may contribute
to the algebraic family Seiberg-Witten invariant
 ${\cal AFSW}_{M_{\delta+1}\times \{t_L\}\mapsto M_{\delta}\times \{t_L\}}(1, 
c_1(L)-2\sum_{1\leq i\leq n}E_i)$ as well. The type $I$ exceptional curves
 begin to contribute to the family invariants when $n\geq 4$. On the other
 hand, when $L$ is $5\delta-1$ very ample, G${\ddot o}$ttsche's argument implies
 implicitly that 
the type $II$ exceptional curves do not contribute to the above family 
invariant ${\cal AFSW}_{M_{\delta+1}\times 
\{t_L\}\mapsto M_{\delta}\times \{t_L\}}(1, 
c_1(L)-2\sum_{1\leq i\leq n}E_i)$. 
But it is a totally different story when $L$ fails to be
 $5\delta-1$ very ample on a ``general'' algebraic surface.

\medskip

(b). From the point of view of family Gromov-Taubes theory, the
 appearance of irreducible type $II$ exceptional curves representing $e_{II}$ 
 pairing negatively with $C-{\bf M}(E)E$ will 
force  effective representatives of $C-{\bf M}(E)E$ lying over
 the existence locus of irreducible effective $e_{II}$ to
 break off a certain multiple of curves representing 
$e_{II}$ and the class $C-{\bf M}(E)E$ can be written as
 $(C-{\bf M}(E)E-e_{II})+(e_{II})$ formally. Potentially they can
 contributes to the family invariant ${\cal AFSW}_{M_{\delta+1}\times 
\{t_L\}\mapsto
 M_{\delta}\times \{t_L\}}(1, C-{\bf M}(E)E)$
as well. 

\medskip

(c). The canonical algebraic obstruction\footnote{Consult [Liu5], [Liu6] for
the detail construction.} bundle 
$\pi_X^{\ast}{\bf W}_{canon}\otimes {\bf H}$ over $X={\bf P}({\bf V}_{canon})$
 was introduced to detect the
 curves with singularities in a complete linear system $|L|$. 
While the base space $M_{\delta}$ parametrizes all the
 possible configurations of singular points and the canonical algebraic 
obstruction
 bundle detects the curves with local multiplicities greater than one, 
any non-reduced curve in 
 the linear system $|L|$ will be detected as singular curves with
 infinitely many singularities. Even if we had subtracted the excess 
contributions from the type $I$ exceptional curves (which are universal), 
the modified family invariant (calculable to be a universal polynomial of
 $C^2$, $C\cdot c_1(K_M), c_1(K_M)^2$ and $c_2(M)$) would still have contained
 the excess contributions from the various multiple covering (corresponding
 to non-reduced curves in algebraic geometry) 
of exceptional curves. 
They are exactly the type $II$
 exceptional curves discussed in definition \ref{defin; two}.

On the other hand, asserting that multiple coverings of type $II$ curves ``may''
contribute potentially to the family invariant 
does not mean that their contributions
 are always nonzero in all possible situations. Our universality theorem
gives us an effective bound universally on the nodal curve counting of 
$L$ such that the type $II$ curves do not contribute to the family
 invariant of $L\otimes {\cal O}(-{\bf M}(E)E)$. In these cases, the
 type $II$ exceptional curve may occur in the possible decompositions
 of the curves dual to $C-{\bf M}(E)E$, but they do not contribute to
the family invariant. This is because these types of decompositions
 involving type $II$ exceptional classes may have a lower 
expected family dimension than the original class $C-{\bf M}(E)E$.

\subsubsection{\bf A Short Discussion of the Vanishing Argument of
 type II Contribution on $K3$}\label{subsubsection; vanishing}

\bigskip

 In the original surface Riemann-Roch
 problem, by Serre duality
 the index of specialty $h^2(L)=h^0(L^{\ast}\otimes K_M)$.
 Then $h^2(L)$ vanishes for sufficiently positive $L$. This corresponds
 to our statement that when $L$ is $5\delta-1$ very ample, type $II$
 curves do not contribute to the family invariant.
 When $K_M$ is trivial (this happens when $M$=$T^4$ or $K3$), it implies the
 vanishing of $h^2(L)$ for any effective $L$. 

 Likewise there is a similar vanishing result for the type $II$ contributions
 to the family invariants when $M=T^4$ or $K3$. In the following, we outline 
the main cause for the vanishing result in a less technical way. 
More details will appear in the paper [Liu7].

\bigskip

 It is well known in Seiberg-Witten theory that $SW({\cal L})=0$ for nontrivial
 $spin_c$ structures ${\cal L}$ over
 all symplectic four-manifolds with trivial canonical bundles. From
 a differential geometric point of view, this vanishing statement
 for $T^4$ and $K3$ are due to
 the existence of Ricci flat metrics on $T^4$ or $K3$ and
 the vanishing result of Witten [W] back to the early days of Seiberg-Witten
 theory. On the
 other hand, in the following we point out directly 
from the point of view of algebraic
 geometry the cause for such a vanishing result.

 Let $M$ be an algebraic $K3$. Let $D$ be an effective
 divisor on $M$ in the linear system $|L|$. 
Then there is the following standard short exact sequence on $M$,

$$0\mapsto {\cal O}_M\mapsto {\cal O}_M(D)\mapsto {\cal O}_D(D)\mapsto 0.$$

If we replace $D$ by the universal divisor ${\bf D}\subset M\times |L|$ and
 push-forward the globalized exact sequence 
$$0\mapsto {\cal O}_{M\times |L|}\mapsto \pi_{\bf D}^{\ast}
{\cal H}\otimes {\cal O}_{M\times |L|}({\bf D})\mapsto \pi_{\bf D}^{\ast}
{\cal H}\otimes {\cal O}_{\bf D}(
{\bf D})\mapsto 0$$

along $\pi_{|L|}:M\times |L|\mapsto |L|$, then we will get the following
 right derived sequence after some simple base change argument,

$$0\mapsto {\bf C}_{|L|}
\mapsto {\bf H}\otimes H^0(M, L)\otimes {\bf C}_{|L|}\mapsto
 {\bf H}\otimes 
{\bf R}^0\pi_{|L|\ast}
{\cal O}_{\bf D}({\bf D})\mapsto \otimes H^1(M, {\cal O}_M)\otimes {\bf C}_{|L|}
\mapsto {\bf H}\otimes H^0(M, L)\otimes {\bf C}_{|L|}$$
$$\mapsto {\bf H}\otimes 
{\bf R}^1\pi_{|L|\ast}\bigl({\cal O}_{\bf D}({\bf D}))
\mapsto {\bf R}^2\pi_{|L|\ast}\bigl({\cal O}_M)\bigr)
\mapsto {\bf H}\otimes {\bf R}^2\pi_{\ast}\bigl({\cal O}({\bf D})\bigr)=0.$$

 This exact sequence has been analyzed in [Liu3] in great detail.
 The key observation here is that by relative Serre duality
 the second derived image 
${\bf R}^2\pi_{\ast}\bigl({\cal O}_M\bigr)
\cong ({\bf R}^0\pi_{\ast}\bigl({\cal O}_M\bigr))^{\ast}$
is isomorphic to 
the trivial line bundle over $|L|$
,${\bf C}_{|L|}$, when
 $M$ has a trivial canonical bundle $K_M$.

 We notice that ${\bf R}^1\pi_{|L|\ast}\bigl({\cal O}_{\bf D}({\bf D}))
\mapsto {\bf R}^2\pi_{|L|\ast}\bigl({\cal O}_M)\bigr)$ is a surjection.
 Consider a Kuranishi-model of ${\bf D}$ with the obstruction bundle 
 ${\bf W}_{obs}$. The moduli space of curves is the zero locus of 
 a defining section of ${\bf W}_{obs}$.
Then by the defining property of the Kuranishi model, 
the fiber of ${\bf W}_{obs}$ above a curve $D\in |L|$ maps surjectively
 onto $H^1(D, {\cal O}_D(D))$.
   
 This implies that the obstruction bundle ${\bf W}_{obs}$ 
of the class $D$ will map surjectively onto the trivial line bundle 
${\bf C}_{|L|}$,
 ${\bf W}_{obs}\mapsto {\bf C}_{|L|}\mapsto 0$. By [F] example 12.1.8 on
 page 216-217, this 
implies that the top Chern class of ${\bf W}_{obs}$ vanishes. 
The fundamental cycle class of the moduli space defined by
 the localized top Chern class of ${\bf W}_{obs}$ along the moduli space of
 curves $\cong |L|$ is zero. In particular, 
the Seiberg-Witten invariant of $D$ is zero.

 This pathetic symptom is exactly why the algebraic (family) Seiberg-Witten
 invariant has been
 introduced [Liu3],
 which removes essentially the troublesome ${\bf C}$ (more generally
 the ${\bf C}^{p_g}$ for $p_g>0$) factor in the
 obstruction bundle and set the obstruction bundle for algebraic Seiberg-Witten
 invariants to be $Ker({\bf W}_{obs}\mapsto {\bf C}_{|L|})$. 

A more intuitive approach is to consider the hyperkahler
twistor family of $K3$. 

 It is well known that $K3$ admits hyperkahler Riemannian 
metrics. I.e. there exists integrable complex structure $I, J, K\in
 End({\bf TM})$ satisfying $I^2=J^2=K^2=-Id_{TM}$ and $IJ+JI=JK+KJ=KI+IK=0$. 
 Then for ${\bf x}=(x, y, z)\in {\bf R}^3, x^2+y^2+z^2=1$, we have an ${\bf S}^2$ 
 family of complex structures $I_{\bf x}=xI+yJ+zK\in End(TM)$.  
 A Ricci-flat metric $g$ is related to the ${\bf S}^2$ family of
 Kahler forms $\omega_{\bf x}$ by 
$g(I_{\bf x}\cdot, \cdot)=\omega_{\bf x}(\cdot, \cdot)$.
 
Consider a special ${\bf S}^2$ hyperkahler family of complex structures on
 the K3 such that $C\in H^2(M, {\bf Z})$ becomes a $(1, 1)$ class
 at $b\in {\bf S}^2$. The class $C$ fails to be of $(1, 1)$ type over 
${\bf S}^2-\{b\}$. So the ``family moduli space'' of curves dual to $C$, 
${\cal M}_C$
 is confined above $b\in {\bf S}^2$.

The tangent space of the thickened base
 ${\bf T}_b{\bf S}^2$ 
maps injectively into the obstruction bundle such that the following
 diagram is commutative,

\[
\begin{array}{ccc}
 {\bf T}_b{\bf S}^2 & \longrightarrow & {\bf W}_{obs}\\
 & \searrow & \Big\downarrow  \\
 & &  {\bf C}_{|L|}
\end{array}
\]

 For the ${\bf S}^2-$thickened hyperkahler family, the obstruction
 bundle ${\bf W}_{obs}/Im{\bf T}_b{\bf S}^2$ of the ${\bf S}^2$ family
 does not have a trivial ${\bf C}$ factor.
 There is no wonder that non-trivial family invariants can be defined for
 such families.

  Let $e_{II}$ be a type $II$ exceptional class of the universal family 
$M_{\delta+1}\mapsto M_{\delta}$. Then the
 family moduli space ${\cal M}_{e_{II}}$ 
over $M_{\delta}$ has the structure of a projectified cone
 over $M_{\delta}\times T(M)$. 
Over the family moduli space ${\cal M}_{e_{II}}$, the
 universal curve (divisor) is denoted by ${\cal D}_{e_{II}}$. 

 Suppose that $C=c_1(L)$ with $C\cdot e_{II}<0$ 
is the class to enumerate the family invariant
 ${\cal AFSW}_{M_{\delta+1}\mapsto M_{\delta}}(1, C-2\sum_i E_i)$.
 Let ${\cal M}_C\mapsto M_{\delta}$ be the family moduli space of curve dual to
 $C$. As a model example 
we investigate the curves dual to $C$
 which are decomposed into components dual to $C-e_{II}$ and $e_{II}$.
 Such curves form a sub-moduli space 
${\cal M}={\cal M}_{C-e_{II}}\times_{M_{\delta}}
{\cal M}_{e_{II}}\subset {\cal M}_C$. Over the locus 
${\cal M}_{C-e_{II}}\times_{M_{\delta}}
{\cal M}_{e_{II}}$, there are a pair of universal curves, 
 ${\cal D}_{C-e_{II}}$ and ${\cal D}_{e_{II}}$ for $C-e_{II}$ and 
 $e_{II}$, respectively. Their sum ${\cal D}_{C-e_{II}}+{\cal D}_{e_{II}}$
 is nothing but the universal curve ${\cal D}_C$ 
of $C$, restricted to the sub-locus 
${\cal M}\subset {\cal M}_C$.

  The key observation is that 
${\bf R}^1\pi_{\ast}\bigl({\cal O}_{{\cal D}_{C-e_{II}}+{\cal D}_{e_{II}}}(
{\cal D}_{C-e_{II}}+{\cal D}_{e_{II}})\bigr)$ is isomorphic to the direct sum 
of a 
rank two sub-sheaf ${\cal O}_{\cal M}\oplus {\cal O}_{\cal M}$ 
 and a coherent sheaf, because
 of the curve decomposition ${\cal D}_C\mapsto {\cal D}_{C-e_{II}}+
 {\cal D}_{e_{II}}$.

Firstly both ${\bf R}^1\pi_{\ast}\bigl({\cal O}_{{\cal D}_{C-e_{II}}}(
{\cal D}_{C-e_{II}}))$ and 
${\bf R}^1\pi_{\ast}\bigl({\cal O}_{{\cal D}_{e_{II}}}({\cal D}_{e_{II}})\bigr)$
 map surjectively onto ${\cal O}_{\cal M}$, following the same discussion 
 earlier.
 Then the mappings \footnote{Both are induced by 
 ${\cal O}_A(A)\mapsto {\cal O}_{A+B}(A+B)$ with $A={\cal D}_{C-e_{II}}$ or
 $A={\cal D}_{e_{II}}$.}
 $${\bf R}^1\pi_{\ast}\bigl({\cal O}_{{\cal D}_{C-e_{II}}}(
{\cal D}_{C-e_{II}})\bigr)\rightarrow 
{\bf R}^1\pi_{\ast}\bigl({\cal O}_{{\cal D}_{C-e_{II}}+{\cal D}_{e_{II}}}(
{\cal D}_{C-e_{II}}+{\cal D}_{e_{II}})\bigr)\leftarrow
{\bf R}^1\pi_{\ast}\bigl({\cal O}_{{\cal D}_{e_{II}}}({\cal D}_{e_{II}})\bigr),$$

send their quotients ${\cal O}_{\cal M}$ injectively into 
 ${\bf R}^1\pi_{\ast}\bigl({\cal O}_{{\cal D}_{C-e_{II}}+{\cal D}_{e_{II}}}(
{\cal D}_{C-e_{II}}+{\cal D}_{e_{II}})\bigr)$. It is not hard to 
check that the direct sum morphism 
${\cal O}_{\cal M}\oplus {\cal O}_{\cal M}\mapsto 
{\bf R}^1\pi_{\ast}\bigl({\cal O}_{{\cal D}_{C-e_{II}}+{\cal D}_{e_{II}}}(
{\cal D}_{C-e_{II}}+{\cal D}_{e_{II}})\bigr)$ is also an injection and the
 corresponding short exact sequence splits. 

 Similar to our earlier discussion, over
the whole family moduli space ${\cal M}_C$ the family obstruction
 bundle ${\bf W}_{obs}$ has a trivial quotient $\cong {\bf C}$. This trivial
 factor is either removed in defining algebraic family Seiberg-Witten
 invariant, or killed by the injection of ${\bf T}_b{\bf S}^2$ if we adopt
 hyperkahler families. When the
 family obstruction bundle ${\bf W}_{obs}$ is restricted to the sub-moduli space
${\cal M}_{C-e_{II}}\times_{M_{\delta}}{\cal M}_{e_{II}}$, it has a 
${\bf C}^2$ quotient 
\footnote{Instead of only one ${\bf C}$.} because of the
 presence of the ${\cal O}_{\cal M}\oplus {\cal O}_{\cal M}$ factor in 
${\bf R}^1\pi_{\ast}\bigl({\cal O}_{{\cal D}_{C-e_{II}}+{\cal D}_{e_{II}}}(
{\cal D}_{C-e_{II}}+{\cal D}_{e_{II}})\bigr)$.

  This implies that the localized contribution of this sub-moduli space
${\cal M}_{C-e_{II}}\times_{M_n}{\cal M}_{e_{II}}$ to the family 
 invariant vanishes, due to
 the presence of additional trivial quotient in the obstruction
 bundle!

 The same argument can be carried out for decompositions with more
 than one type $II$ components and we find that
 the type II curves contribute
 trivially to the family invariants of $K3$. This is the reason why the
 type $II$ correction terms never appear in the nodal curve counting on $K3$.

\medskip

\bigskip

\section{The 
Outline of Some Key Ideas in the Algebraic Proof of
 Universality Theorem}\label{section; outline}

\bigskip

 Even though the algebraic 
proof of universality theorem in [Liu6] has been reduced to 
a purely algebraic argument and hence the conceptual dependence on
 the symplectic $''SW=Gr''$ used in [Liu1] has been relieved, 
the argument is still of a considerable length.
 To guide the reader through the paper, the goal of the section is to provide
 a short sketch of the ideas involved in the proof, the various difficulties it
 has overcome, etc., in less technical terms.

 Let $M_{\delta}$ denote the $n-$th universal space associated to the algebraic
surface $M$. Then $f_{\delta}: M_{\delta+1}\mapsto M_{\delta}$ is smooth of relative dimension
 two and is the universal space parametrizing the $n-$consecutive codimension
 two blowing ups of
 $M$.
 In [Liu3] section 5.1 we have introduced the concept of canonical algebraic family
Kuranishi model\footnote{ Here $C=c_1(L)$. We use $C$ in the paper to 
allow the extension to non-linear systems instead of linear system $|L|$.}
 of $C-{\bf M}(E)E$. Assuming that ${\cal R}^i\pi_{\ast}\bigl({\cal E}_C\bigr)=0$
for $i>0$, the pair of vector bundles 
${\bf V}_{canon}, {\bf W}_{canon}$ were introduced,
 where ${\cal V}_{canon}={\cal R}^0\pi_{\ast}\bigl({\cal E}_C\bigr)$ and
 ${\cal W}_{canon}={\cal R}^0\pi_{\ast}\bigl({\cal O}_{{\bf M}(E)E}\otimes
 {\cal E}_C\bigr)$ are the corresponding locally free sheaves of sections.
 The sheaf morphism 
 ${\cal R}^0\pi_{\ast}\bigl({\cal E}_C\bigr)\mapsto 
{\cal R}^0\pi_{\ast}\bigl({\cal O}_{{\bf M}(E)E}\otimes
 {\cal E}_C\bigr)$ between the locally free sheaves (see definition 5.3 of
 [Liu3] for 
more details) induces a bundle map between the vector bundles
 ${\bf V}_{canon}\mapsto {\bf W}_{canon}$. The vector bundle 
${\bf V}_{canon}\mapsto M_{\delta}\times T(M)$
parametrizes the non-linear system of curves dual to $C$ in the family
 $M_{\delta+1}\times T(M)\mapsto M_{\delta}\times T(M)$. A curve (corresponding to
 a ray in a fiber of ${\bf V}_{canon}$) dual to $C$ can be splitted into
 components dual to
 $C-{\bf M}(E)E$ and dual to ${\bf M}(E)E$ if and only if the 
 restriction map to ${\bf M}(E)E$ of the curve dual to $C$ vanishes along the 
whole ${\bf M}(E)E$, i.e. the ray corresponding to the definition section
 of the curve is in the kernel $Ker({\bf V}_{canon}\mapsto
 {\bf W}_{canon})$ of the bundle map\footnote{This canonical bundle map is
 constructed by the restriction morphism to ${\bf M}(E)E$.}
 ${\bf V}_{canon}\mapsto {\bf W}_{canon}$.

 The bundle map ${\bf V}_{canon}\mapsto {\bf W}_{canon}$ induces a bundle
 map ${\bf H}^{\ast}\mapsto \pi_X^{\ast}{\bf W}_{canon}$ over 
$X={\bf P}({\bf V}_{canon})$, which is nothing but a global section of
 $\pi_X^{\ast}{\bf W}_{canon}\otimes {\bf H}$.
  Thus the family moduli space of $C-{\bf M}(E)E$ can be thought to be
a projectified cone embedded in $X={\bf P}({\bf V}_{canon})$, defined
 by a canonical section $s_{canon}\in \Gamma(X, \pi_X^{\ast}{\bf W}_{canon}\otimes
 {\bf H})$ induced by ${\bf V}_{canon}\mapsto {\bf W}_{canon}$. It is
 the union of the projectified (non-)linear systems of the fibers. 
 Thus, we may rewrite ${\cal M}_{C-{\bf M}(E)E}=Z(s_{canon})$.
  If we want to count curves in the linear system $|L|$ instead of 
the non-linear system, then we work with 
${\cal M}_{C-{\bf M}(E)E}\times_{T(M)}\{t_L\}=Z(s_{canon})\times_{T(M)}\{t_L\}$.
 Over here $t_L\in T(M)$ is a point in the component of 
Picard variety $T(M)\cong Pic^0(M)$ corresponding to $L$.

  By intersection theory [F] construction,
 the localized top Chern class defines a cycle class
 supported in $Z(s_{canon})$, it is 
denoted by ${\bf Z}(s_{canon})$, following
 the notations of [F].
 Consider the top open stratum $Y_{\gamma_{\delta}}$ of the universal space
 $M_{\delta}$ which parametrizes the
 ordered distinct $n$ points on $M$. Then $M_{\delta}-Y(\gamma_{\delta})$ is a closed
sub-scheme of $M_{\delta}$. In the following, we denote 
$Y=X\times_{M_{\delta}}(M_{\delta}-Y_{\gamma_{\delta}})$.

 We can split $Z(s_{canon})$ into the closed 
 $Z_o=\overline{Z(s_{canon})\times_{M_{\delta}}Y_{\gamma_{\delta}}}$ and 
 the closed $Z_i=\overline{Z(s_{canon})-Z_o}$

\noindent 
$\overline{Z(s_{canon})-\overline{Z(s_{canon})\times_{M_{\delta}}Y_{\gamma_{\delta}}}}
\subset Z(s_{canon})\times_{M_{\delta}}Y$. The latter is the union of all the irreducible
 components in $Z(s_{canon})\times_{M_{\delta}}Y$.

 We also know that $X-Y=X\times_{M_{\delta}}Y_{\gamma_{\delta}}\mapsto 
 X$ is the
 top open stratum of $X$, where the stratification on $X$ has been
 induced from the admissible stratification \footnote{See section 2 of
 [Liu6] for more
 details about how to stratify $M_{\delta}$ by $Y_{\Gamma}$, $\Gamma\in adm(n)$.} 
of $M_{\delta}$ by the surjection 
$X\mapsto M_{\delta}$. On the other hand, the standard exact sequence 
 for $U=Z(s_{canon})-Z(s_{canon})\times_{M_{\delta}}Y$, 

$${\cal Z}_{\cdot}(Y\cap Z(s_{canon}))
\mapsto {\cal Z}_{\cdot}(Z(s_{canon}))\stackrel{j^{\ast}}{\mapsto}
 {\cal Z}_{\cdot}(Z(s_{canon})-Z(s_{canon})\cap Y)
\mapsto 0$$

and the induced exact sequence on their quotients ${\cal A}_{\cdot}$,

$${\cal A}_{\cdot}(Y\cap Z(s_{canon}))
\mapsto {\cal A}_{\cdot}(Z(s_{canon}))\stackrel{j^{\ast}}{\mapsto}
 {\cal A}_{\cdot}(Z(s_{canon})-Z(s_{canon})\cap Y)
\mapsto 0$$

 on page 21 of [F] imply that ${\bf Z}(s_{canon})\in 
{\cal A}_{\cdot}(Z(s_{canon}))$ can be splitted into a unique
 component in ${\cal A}_{\cdot}(Z(s_{canon})\cap Y)$ and 
 a unique image in ${\cal A}_{\cdot}(Z(s_{canon})-Z(s_{canon})\cap Y)$.

  As $Z_o$ and $Z_i$ are both unions of irreducible components of 
 $Z(s_{canon})$,
 the decomposition $Z(s_{canon})=Z_o\cup Z_i$ induces a 
decomposition on the normal cones ${\bf C}_{Z(s_{canon})}X={\bf C}_{Z_o}X\cup
{\bf C}_{Z_i}X$, which induces a corresponding decomposition of
 $s({\bf C}_{Z(s_{canon})}X)\in {\cal A}_{\cdot}(Z(s_{canon}))$ into
 $s({\bf C}_{Z_o}X)\in {\cal A}_{\cdot}(Z_o)$ and 
 $s({\bf C}_{Z_i}X)\in {\cal A}_{\cdot}(Z_i)$. By capping with
 $c_{total}({\bf H}\otimes \pi_X^{\ast}W_{canon}|_{Z(s_{canon})})$, we get a
 canonical decomposition of ${\bf Z}(s_{canon})$ into a unique component
 in ${\cal A}_{\cdot}(Z(s_{canon})\cap Y)$, called localized contribution
 of top Chern class ${\bf Z}_{Z_o}(s_{canon})$ and
 a component extended from ${\cal A}(Z(s_{canon})-Z(s_{canon})\cap Y)$.

\begin{defin}\label{defin; extension}
Let $\eta\in {\cal A}_{\cdot}(Z(s_{canon})-Z(s_{canon})\cap Y)=
{\cal A}_{\cdot}(Z_o-Z_o\cap Y)$. Lift $\eta$ to an explicit closed cycle
 in ${\cal Z}_{\cdot}(Z_o-Z_o\cap Y)$. Then such a lifting can be
 extended to ${\cal Z}_{\cdot}(Z(s_{canon}))$. 
Define the resulting extension of $\eta$ to ${\cal A}_{\cdot}(Z(s_{canon}))$ by
 $\overline{\eta}$.
\end{defin}

The extension $\overline{\eta}$ a priori depends on the lifting of $\eta$ from 
the cycle class group to
 ${\cal Z}_{\cdot}(Z_o-Z_o\cap Y)$. On the other hand, for 
$\eta=j^{\ast}{\bf Z}(s_{canon})$ we can take $\overline{\eta}$ to be
 ${\bf Z}_{Z_o}(s_{canon})$ $=\{c_{total}(\pi_X^{\ast}{\bf W}_{canon}\otimes
 {\bf H}|_{Z_o})$$\cap s_{total}(Z_o, X)
\}_{dim_{\bf C}X-rank_{\bf C}{\bf W}_{canon}}$ and it is unique.

Suppose that we had had the equality 
$\overline{j^{\ast}{\bf Z}(s_{canon})}={\bf Z}(s_{canon})$, then
 ${\bf Z}(s_{canon})$ would have been viewed as an extension of the
 cycle class $j^{\ast}{\bf Z}(s_{canon})$ into $Z(s_{canon})$ 
(by taking the closure!) 
and it would have represented the fundamental cycle class of the ``moduli space of
 curves from the fiberwise sections of 
${\cal E}_C\otimes {\cal O}(-{\bf M}(E)E)$. Each curve in
 $Z(s_{canon})\times_{M_{\delta}}Y_{\gamma_{\delta}}$ projects into 
 $M$ by
 the restricted 
blowing down map $M_{\delta+1}\times_{M_{\delta}}Y_{\gamma_{\delta}}\mapsto M\times Y_{\gamma_{\delta}}$
 and produces  
 curves with at least $n$ singularities.  One may
 consider the intersection number by capping 
the fundamental class ${\bf Z}(s_{canon})$ with 
 $c_1({\bf H})^{p_g-q+{(C-{\bf M}(E)E)^2-(C-{\bf M}(E)E)\cdot c_1(K_{M_{\delta+1}/M_{\delta}})
\over 2}}$
 and the answer can be shown to be a universal degree $n$ polynomial of
 $c_1(L)^2, c_1(L)\cdot c_1(M), c_1^2(M)$ and $c_2(M)$, by applying 
the family blowup formula [Liu3] inductively.

 The possible failure of the equality
${\bf Z}(s_{canon})=\overline{j^{\ast}{\bf Z}(s_{canon})}$
 indicates that there are excess cycle classes 
${\bf Z}(s_{canon})-\overline{j^{\ast}{\bf Z}(s_{canon})}$ localized in
 ${\cal A}_{\cdot}(Y\cap Z(s_{canon}))$ 
which also contributes to ${\bf Z}(s_{canon})$ besides
 the generic component 
$j^{\ast}{\bf Z}(s_{canon})$. Thus it becomes a subtle issue to separate
 $\overline{j^{\ast}{\bf Z}(s_{canon})}$ from the localized top Chern class 
${\bf Z}(s_{canon})$, which
 can be pushed-forward into $X$ as the global object--the top Chern class 
$c_{top}(\pi_X^{\ast}{\bf W}_{canon}\otimes {\bf H})$. 

 Different groups of people had chosen quite different paths to either 
address or to bypass
 this issue. In [Got] G${\ddot o}$ttsche had proposed a different approach using 
 Hilbert schemes of surfaces. Instead of using universal spaces and counting
 the resolved curves, he works with the singular curves themselves instead of
 the resolved smooth curves. 
When $L$ is $5\delta-1$ very ample, he could construct
 a finite scheme (see [Got]) 
in an ambient scheme birational to some Hilbert scheme of $M$, 
representing the fundamental cycle class of moduli space of nodal curves in generic
 $n$ dimensional linear sub-system. In his approach,
 he completely avoided the above problem. On the other hand,
 the topological nature of the resulting ``number of nodal curves''
 has been less transparent.

 On the other hand, Vainsencher [V], Kleiman-Piene 
had adopted the universal space approach
 without using the admissible stratification
explicitly
and fought with the problem directly.
 By assuming that $L$ is sufficiently
 very ample and $\delta\leq 6$, $\leq 8$, respectively, they were able to 
 prove regularity results on the $Z(s_{canon})$ and identify
 the non-trivial excess contribution of intersection numbers from 
${\cal A}_{\cdot}(Z(s_{canon})\cap Y)$
 in the forms of sums of excess contributions of singular curves from other 
types of
 non-nodal singularities\footnote{See e.g. [V] for more details and
 examples.}.
 
 \bigskip

 In [Liu1], [Liu6] the author has taken a new approach to tackle the problem.
 We realize that the problem is in fact the natural prolongation of 
surface Riemann-Roch formula. Then we migrate tools from differential
 topology and symplectic geometry to algebraic geometry in resolving the
 problem. The family blowup formula and the family switching formula
 provide the necessary bridges between our non-linear enumeration problem
 and the Grothendieck Riemann-Roch theorem.

\subsection{\bf Some Problems Addressed in Our Paper} \label{subsection; address}

\bigskip

 In the following we list all the key problems we address in the paper [Liu6],

\medskip

\noindent 
{\bf Problem 1}: As has been commented earlier, we do not expect that
for all $L$ and the general $M$, 
the numbers of nodal curves in linear sub-system of $|L|$ 
 are topological numbers. What are the effective range of $L$ in which the
 numbers of nodal curves become topological?

If one always has to raise $L$ to an extremely high power of a
very ample line bundle, without an effective
 control over its very-ampleness, 
it limits the application of the universal formulae to enumerative 
predictions. 

\noindent 
{\bf Problem 2}: It is rather unclear from a geometric prospective that the
 excess 
localized contributions of the cycle classes in ${\cal A}_{\cdot}(Z(s_{canon})\cap 
Y)$ are ``intrinsic''. I.e. it is not clear at all that 
the cycle classes localized inside ${\cal A}_{\cdot}(Y)$ 
 are ``invariant'' to the deformations of the complex structures of $M$ 
 and/or the holomorphic structures of $L$. Not to mention that it may or
 may not have
any intrinsic properties across different algebraic surfaces.

 Consider an ideal situation that $s_{canon}$ can be slightly deformed
 (in the algebraic category) into a new cross section $s'$ with the zero locus
 $Z(s')$. It is usually not true that ${\bf Z}(s')-\overline{j^{\ast}{\bf Z}(s')}$
 is equal to ${\bf Z}(s_{canon})-\overline{j^{\ast}{\bf Z}(s_{canon})}$.

 In what sense does the excess localized contributions to top Chern class
 ${\bf Z}(s_{canon})-\overline{j^{\ast}{\bf Z}(s_{canon})}$ have any 
intrinsic geometric meaning?

\bigskip

 On the other hand,

\noindent {\bf Problem 3}: The space $Y=X\times_{M_{\delta}}(M_{\delta}-Y_{\gamma})$ 
itself is not smooth. Potentially this
adds to the technical issues in identifying the excess cycle class localized
 in $Z(s_{canon})\cap Y$. In general the locus $Z_i\subset Y$ is highly
 singular. How do we study enumerative geometry upon such singular objects?

\bigskip

 Our strategy to simplify {\bf Problem 3} 
is to use the admissible stratification on the
universal space $M_{\delta}=\coprod_{\Gamma\in adm(n)}Y_{\Gamma}$ and rewrite $Y$ as
 $X\times_{M_{\delta}}(\cup_{\Gamma\in adm(n); \Gamma\not=\gamma_{\delta}}Y(\Gamma))$,
 which is a finite union of smooth subspaces $X\times_{M_{\delta}}Y(\Gamma)$ in $X$. 

 In this way we have replaced a single non-smooth $Y$ by a hierachy of
 smooth $X\times_{M_{\delta}}Y(\Gamma)$ and different $X\times_{M_{\delta}}Y(\Gamma)$
 may intersect each other.

\bigskip 

One key observation in the algebraic proof of universality theorem 
is the following simple fact: 
The algebraic analogue of the Kuranishi model type perturbation argument in
 the differentiable or symplectic 
category is exactly the localized top Chern class of the zero
 locus of a section
 of the algebraic obstruction vector bundle. 
We generalize this concept slightly and call the following expression
 $\{c_{total}(E)\cap s_{total}(Z, X)\cap [X]\}_{dim_{\bf C}X-rank_{\bf C}E}$
 the localized contribution of top Chern class for $Z\subset Z(s)$ and some
 $s\in \Gamma(X, E)$.

 On the other hand, one may blow up $X$ along $Z$ with an exceptional divisor
 $D$. Then the above expression of localized contribution of top Chern
class can be viewed as the push-forward \footnote{by using the
definition of total Segre class of the normal cone ${\bf C}_ZX$.}
 of $\sum_{1\leq i\leq 
 rank_{\bf C}E} c_{rank_{\bf C}E-i}(E)\cap (-1)^{i-1}
D^{i-1}[D]$, which appears naturally in the residual intersection formula
 of top Chern class. See example 14.1.4 on page 245 of [F] for more details.

 Our algebraic geometric 
construction to enumerate the localized contribution of top Chern classes 
along $Z(s_{canon})\cap Y$ 
is to blow up $X$  repeatedly along the sub-loci.
 But we do not blow up $X$ along $Z(s_{canon})\cap Y$
 at once because it leads to an un-identifiable
 localized contribution of top
Chern classes and is not helpful to us. Instead 
our scheme to tackle this problem requires us to rewrite
 $Z(s_{canon})\cap Y$ as 
a seemingly more complicated union 
$\cup_{\Gamma\in \Delta(n)-\{\gamma_{\delta}\}}Z(s_{canon})\times_{M_{\delta}}Y(\Gamma)$.
 Notice that quite a lot of admissible strata $Y(\Gamma)$ with $\Gamma\in adm(n)$
 \footnote{The set $adm(n)$ is the set of $n$-vertex admissible graphs. 
 See [Liu6] section 2 for its definition.}
 
have been thrown away from the union as they are in the closure of
 the other admissible strata. One can show that we only need to consider those
 $\Gamma\in \Delta(n)$, satisfying special maximality conditions\footnote{Consult
 definition 6 on page 47 of the paper [Liu6].}.
The set $\Delta(n)$ collects the admissible graphs such that their associated
admissible strata are kept in the union. Then we blow up $X$ inductively along 
the various $Z(s_{canon})\times_{M_{\delta}}Y(\Gamma)$, $\Gamma\in\Delta(n)$, and apply
 the residual intersection formula of top Chern class to 
$\pi_X^{\ast}{\bf W}_{canon}\otimes {\bf H}$ and to the canonical
 section $s_{canon}$ defining ${\cal M}_{C-{\bf M}(E)E}$.

\medskip

  But this approach also rises three additional 
issues that we have to address, in order to identify
 all these localized contributions of top Chern classes.

\medskip

\noindent {\bf Problem 4}: Is the contribution of localized top Chern classes 
along $Z(s_{canon})\times_{M_{\delta}}Y(\Gamma)$ by an inductive 
application of residual intersection theory of top Chern classes 
sensitive to the ordering
 that we choose for the blowing ups of the various 
$Z(s_{canon})\times_{M_{\delta}}Y(\Gamma')$ applied prior to the given 
 $Z(s_{canon})\times_{M_{\delta}}Y(\Gamma)$?

 and,

\medskip

\noindent {\bf Problem 5}: The canonical algebraic obstruction
 bundle $\pi_X^{\ast}{\bf W}_{canon}\otimes {\bf H}$ defining algebraic
 family Seiberg-Witten invariants of $C-{\bf M}(E)E$ 
 gets modified into $\pi_X^{\ast}{\bf W}_{canon}\otimes {\bf H}\otimes 
 {\cal O}(-D)$ even after a single usage of residual intersection formula of 
top Chern class. The symbol here $D$ stands for the exceptional divisor of
 the blowing up.

 Suppose we blow up along each $Z(s_{canon})\times_{M_{\delta}}Y(\Gamma)$
 individually, potentially we have to perform millions of blowing ups on $X$ when
 $n$ goes large.

\medskip 

 Are we still able to work with the seemingly complicated residual
obstruction vector bundle \footnote{They are denoted as 
 $\pi_X^{\ast}{\bf W}_{canon}\otimes {\bf H}\otimes_{\Gamma'\in I_{\Gamma}}
 {\cal O}(-D_{\Gamma'})$ in the paper [Liu6].} which get modified and identify
 the various localized contributions of top Chern class localized in different
 $Z(s_{canon})\times_{M_{\delta}}Y(\Gamma)$,  $\Gamma\in \Delta(n)-\{\gamma_{\delta}\}$?

\medskip

\noindent {\bf Problem 6}: As we have mentioned that 
different $Y(\Gamma)$, $\Gamma\in \Delta(n)-
\{\gamma_{\delta}\}$, may intersect non-trivially. As a consequence
 different sub-loci $Z(s_{canon})\times_{M_{\delta}}Y(\Gamma)$ can touch.
 There can be a potential danger of ``over counting''. How do 
 we avoid over-counting in our scheme? 

\bigskip

\subsection{\bf Responses to the Above Problems}

\bigskip

 All these six problems have addressed in the paper [Liu6], which is
 responsible of its length and complicated notations.
 Let us sketch how we have addressed these issues in the paper [Liu6],
 phrased in a less technical term. Hopefully it can provide the reader
 a guide to read the long paper.

\medskip

\noindent Response to {\bf Problem 1}: While identifying the cycle class
 $\overline{j^{\ast}{\bf Z}(s_{canon})}$ geometrically, we do not attempt to
 raise $L$ to a very high power to control the regularity of the scheme 
 $Z(s_{canon})\times_{M_{\delta}}Y$.  Even though the regularity of
 $Z(s_{canon})\times_{M_{\delta}}Y$ (or $Z(s_{canon})\times_{M_{\delta}}Y\cap V$)
 \footnote{here $V$ is the generic $\delta$-dimensional linear sub-system}
 may be helpful in identifying the excess contribution
 $({\bf Z}(s_{canon})-\overline{j^{\ast}{\bf Z}(s_{canon})})\cap 
c_1({\bf H})^{rank_{\bf C}{\bf V}_{canon}-rank_{\bf C}{\bf W}_{canon}-1+2n}\cap
 [t_L]$, we construct new machineries to 
bypass the problem--by using the combinatorial structure on the universal
 spaces suggested by Gromog-Taubes theory in symplectic geometry.

 We only use a slightly strengthened form
 of G${\ddot o}$ttsche's argument to control the locus $Z(s_{canon})\times_{M_{\delta}} 
Y_{\gamma_{\delta}}$.
 Under the $5\delta-1$ very ampleness assumption on $L$, 
 $V\cap Z(s_{canon})\times_{M_{\delta}} Y_{\gamma_{\delta}}$ is a finite sub-scheme in $X$
 for generic $\delta$ dimensional 
$V$. The main theme of the paper [Liu6] is to identify
 the excess contributions to the family invariant ${\cal AFSW}_{M_{\delta+1}\times 
\{t_L\}\mapsto M_{\delta}\times \{t_L\}}(1, c_1(L)-2\sum_{1\leq i\leq n}E_i)$
using only the smoothness of the family moduli spaces of type $I$ exceptional
 classes and their intersections: $Y(\Gamma)$.
 We do not try to control the regularity of $Z(s_{canon})\times_{M_{\delta}}Y$.

\medskip

\noindent Response to {\bf Problem 2}: The original problem of enumerating the
 localized contribution of top Chern classes is
 formulated in terms of intersection theory [F]. Yet
 to answer {\bf problem 2}, the ideas inspired from family Gromov-Taubes 
 theory have played crucial roles. The algebraic curves representing the
 class $C-{\bf M}(E)E$ within the family $M_{\delta+1}\times T(M)\mapsto M_{\delta}\times 
 T(M)$ may fail to be irreducible and may break into more
 than one irreducible component. 
 While the curves breaking into different components is a purely geometric
 phenomenon, we have found a topological constraint which forces the
 degeneration to occur.

\begin{lemm}\label{lemm; break}
Let ${\bf e}$ be an irreducible exceptional curve representing the exceptiona
 class $e$ over $b\in M_{\delta}$ in the family $M_{\delta+1}\mapsto M_{\delta}$. 
 Suppose that $(C-{\bf M}(E)E)\cdot e<0$, then any effective
 representative of $C-{\bf M}(E)E$ over $b$ has to contain ${\bf e}$ as
 one of its irreducible components.
\end{lemm}

\noindent Proof:  If the curve $\Sigma$ representing $C-{\bf M}(E)E$ lies in the
 same fiber $M_{\delta+1}\times_{M_{\delta}}\{b\}$ as ${\bf e}$ but does not contain
 ${\bf e}$ as one of its irreducible components. Then all irreducible 
components of $\Sigma$ intersect non-negatively with ${\bf e}$. As a
 consequence their sum $=(C-{\bf M}(E)E)\cdot e\geq 0$, violating the
 assumption! $\Box$ 

Over the various locally closed and smooth admissible strata
 $Y_{\Gamma}$, $\Gamma\in \Delta(n)-\{\gamma_{\delta}\}$, those
type $I$ exceptional classes 
 $e_i$ with $e_i^2<-1$, effective and irreducible over $Y_{\Gamma}$, 
 \footnote{Suppose $j_i$ denote the direct descendent indexes of $i$ 
in $\Gamma$, then $e_i=E_i-\sum_{j_i}E_{j_i}$.} pair negatively
\footnote{They are denoted as $e_{k_i}, 1\leq i\leq p$ in
 [Liu6].} with $C-{\bf M}(E)E$. When $e_i$ is represented 
 by irreducible\footnote{It is the case above points in $Y_{\Gamma}$.}
 ${\bf P}^1$, by lemma \ref{lemm; break} 
the condition $(C-{\bf M}(E)E)\cdot e_i<0$
 forces any algebraic curve representing $C-{\bf M}(E)E$ above the same fiber
 to break off at least an irreducible ${\bf P}^1$ component
 representing $e_i$.
 This suggests
 that the canonical algebraic family obstruction bundles \footnote{Refer to
 proposition 9 of [Liu5] for more details.} 
 $\pi_X^{\ast}{\bf W}_{canon}^{\circ}\otimes {\bf H}$
 of the new class $C-{\bf M}(E)E-\sum_{e_i\cdot (C-{\bf M}(E)E)<0}e_i$
 over $Y(\Gamma)$ will be involved. In [Liu5] and [Liu6], we discuss 
extensively the relationship of $\pi_X^{\ast}{\bf W}_{canon}^{\circ}\otimes
 {\bf H}$, its canonical section $s_{canon}^{\circ}$ and 
$\pi_X^{\ast}{\bf W}_{canon}\otimes
 {\bf H}$, with its canonical section $s_{canon}$. The unique
 factorization of curves from $C-{\bf M}(E)E$ to 
$C-{\bf M}(E)E-\sum_{e_i\cdot (C-{\bf M}(E)E)<0}e_i$ can be formulated
 as the isomorphism $Z(s_{canon}^{\circ})\times_{M_{\delta}}Y_{\Gamma}\cong 
Z(s_{canon})\times_{M_{\delta}}Y_{\Gamma}$ of moduli spaces of curves. 
 The algebraic counter-part of family switching formula [Liu5] 
suggests a natural bundle morphism 

$$\pi_X^{\ast}{\bf W}_{canon}^{\circ}\otimes {\bf H}|_{X\times_{M_{\delta}}Y(\Gamma)}
\longrightarrow 
\pi_X^{\ast}{\bf W}_{canon}\otimes {\bf H}|_{X\times_{M_{\delta}}Y(\Gamma)}$$ whose
 injectivity over $X\times_{M_{\delta}}Y_{\Gamma}$ is responsible for the 
above isomorphism of moduli space of curves.

 The above algebraic structures are suggested by the perturbation
 argument on Kuranishi models in the differentiable category [Liu1].

 We analyze the cycle class
 localized in $Z(s_{canon})\cap Y$ and show that it can be 
 identified with a huge 
sum of terms (one for each $\Gamma\in \Delta(n)-\{\gamma_{\delta}\}$,
 through the residual intersection
 formula of top Chern class) of certain cycle classes associated to the different
 obstruction bundles $\pi_X^{\ast}{\bf W}_{canon}^{\circ}\otimes
 {\bf H}$, one for each statum $Y(\Gamma)$.
 This procedure produces intersection numbers known as
 the modified algebraic family Seiberg-Witten invariants of 
$C-{\bf M}(E)E-\sum_{e_i\cdot (C-{\bf M}(E)E)<0}e_i$,
 denoted as 
${\cal AFSW}_{M_{\delta+1}\times_{M_{\delta}}Y(\Gamma)\mapsto Y(\Gamma)}^{\ast}
(c_{total}(\tau_{\Gamma}), C-{\bf M}(E)E-\sum_{e_i\cdot (C-{\bf M}(E)E)<0}e_i)$.

 Once we can identify the intersection number ``localized'' to each
 $X\times_{M_{\delta}}Y_{\Gamma}$ with the modified invariant,
 the fact that all these modified algebraic family Seiberg-Witten
 invariants are ``topological objects'', i.e. can be expressed as
 characteristic classes implies that all the intersection numbers attached
 to $Z(s_{canon})\times_{M_{\delta}}Y(\Gamma)\cap V$
 are topological objects. 
 
\medskip

\noindent Response to {\bf Problem 3}: Even though $Y$ itself is not smooth, $Y$
 can be identified with the unions of the various 
$X\times_{M_{\delta}}Y(\Gamma)$, which 
are all smooth. The type $I$ exceptional classes effective over $Y(\Gamma)$ have
 played an essential role in answering Problem 2. On the other hand,
 each $Y(\Gamma)$ can be intrepreted as the locus of co-existence of 
 the type $I$ exceptional class $e_i$ defined by $\Gamma$ (see section
 2, proposition 4 of [Liu6]),

 $$Y(\Gamma)=\cap_{1\leq i\leq n}Y(\Gamma_{e_i}).$$

 The above intersection is a regular intersection of smooth spaces
 $Y(\Gamma_{e_i})$. The $\Gamma_{e_i}$ is an admissible sub-graph of
 $\Gamma$ with one-edges from $i-$th vertex to all its direct discendents
 in $\Gamma$. They are called fan-like admissible graphs in [Liu6].
 The fan-like admissible graphs $\Gamma_{e_i}$ and $\Gamma$ give us
 combinatorial tools to express the enumerative geometric data.

 For each $Y(\Gamma)\subset M_{\delta}$
 the normal bundle ${\bf N}_{Y(\Gamma)}M_{\delta}$ is isomorphic to
 the restriction of the direct sum of normal bundles of $Y(\Gamma_{e_i})\subset
 M_{\delta}$, which are isomorphic to the restriction of 
canonical algebraic obstruction bundles
 of $e_i$, $1\leq i\leq n$ to $Y(\Gamma)$.
 This information is crucial for us to identify the individual
 terms of excess contributions from the residual intersection formula.

 The key is the following four-term sheaf exact sequence

$$\hskip -.5in 
0\mapsto 
{\cal R}^0\pi_{\ast}\bigl({\cal O}_{\sum_{1\leq i\leq \Xi_{k_i}}}\otimes
 {\cal E}_{C-{\bf M}(E)E}\bigr)\mapsto {\cal R}^0\pi_{\ast}\bigl(
{\cal O}_{{\bf M}(E)E+\sum_{1\leq i\leq p}\Xi_{k_i}}\otimes {\cal E}_C
\bigr)|_{Y(\Gamma)\times T(M)}\mapsto $$ 
$$\hskip -.7in
{\cal R}^0\pi_{\ast}\bigl({\cal O}_{{\bf M}(E)E}
\otimes {\cal E}_C\bigr)|_{Y(\Gamma)\times T(M)}\mapsto 
{\cal R}^1\pi_{\ast}\bigl({\cal O}_{\sum_{1\leq i\leq \Xi_{k_i}}}\otimes
 {\cal E}_{C-{\bf M}(E)E}\bigr)
\mapsto 0,$$

 where $\Xi_{k_i}$ are the universal type $I$ exceptional curve fibrations
 representing $e_{k_i}$.

The second and the third terms in the above sequence are locally free and
 their associated vector bundles are ${\bf W}_{canon}^{\circ}$
 and ${\bf W}_{canon}$, respectively.

 This four-term exact sequence allows us to ``transfer'' the localized contribution
 of top Chern class of $\pi_X^{\ast}{\bf W}_{canon}\otimes {\bf H}$ 
along $Z(s_{canon})\times_{M_{\delta}}Y(\Gamma)$ to 
 the top Chern class of $\pi_X^{\ast}{\bf W}_{canon}^{\circ}\otimes 
{\bf H}$ and chern classes of some additional obstruction virtual
 bundle $\tau_{\Gamma}$ (see definition 10 of [Liu6]).

\medskip

\noindent Response to {\bf Problem 4}: In proposition 16 of [Liu6], 
we show that suitable permuting
 the blowing up orders or some special collapsing/grouping of 
family moduli spaces does not affect the answer of the evaluation. 
 Let $P$ be an index subset of $\Delta(n)-\{\gamma_{\delta}\}$. 
The key observation is that no matter which 
$Z(s_{canon})\times_{M_{\delta}}Y(\Gamma')$
 we have blown up earlier, when we push-forward the sum of the various localized
contributions of top Chern classes (indexed by elemenets in $P$) 
to $X$, the total sum is always equal
to $\{c_{total}(\pi_X^{\ast}{\bf W}_{canon}\otimes {\bf H})\cap
 s_{total}(Z(s_{canon})\times_{M_{\delta}}(\cup_{\Gamma\in P}Y(\Gamma)), X)
\cap [X]\}_{dim_{\bf C}X-rank_{\bf C}{\bf W}_{canon}}$. 
The rationale behind this statement is based on the
 well known property that total Segre class is invariant under
 proper birational push-forward (see proposition 4.2 on page 74 of [F]).
 
  This observation allows us to permute the blowup orderings ahead of
 any particular blowing up while identifying
 a particular localized contribution of the top Chern class.

\begin{rem}
 The permutation of the blowing up orderings
 depend on each individual $Y(\Gamma)$ 
 we focus upon and is
 not universal. 
\end{rem}

\medskip

\noindent Response to {\bf Problem 5}: {\bf Problem 5} is conceptually the most
 challenging one to answer in our approach. 
Our key observation, proposition 9 of [Liu6], answers the question in a
 slightly surprising way. Quite opposite to the naive intuition,
 the repeated blowing ups does not spoil
 our enumeration program. In fact it is necessary
 in identifying the localized contribution
 of top Chern classes!

In the earlier responses to {\bf problem 2 \& 3}, 
we has made use of ${\bf W}_{canon}^{\circ}$
 to replace ${\bf W}_{canon}$. On the other hand, the bundle map
 $\pi_X^{\ast}
{\bf W}_{canon}^{\circ}\otimes {\bf H}|_{X\times_{M_{\delta}}Y(\Gamma)}\mapsto
\pi_X^{\ast}{\bf W}_{canon}\otimes {\bf H}|_{X\times_{M_{\delta}}Y(\Gamma)}$ 
is generically injective over $X\times_{M_{\delta}}Y_{\Gamma}$ but
 may fail to be injective on a closed subset of 
$X\times_{M_{\delta}}\bigl(Y(\Gamma)-Y_{\Gamma}\bigr)$.
 The failure of the bundle injection is related to the appearance of 
some other $Y(\Gamma')$, $\Gamma\succ \Gamma'$; $Y(\Gamma)\cap Y(\Gamma')\not=
\emptyset$\footnote{Consult corollary 3 on page 40 of [Liu5]
 for more details.}. The failure of the bundle map to be injective 
throughout $X\times_{M_{\delta}} Y(\Gamma)$ forbids us to
 replace $c_{top}(\pi_X^{\ast}{\bf W}_{canon}|_{Y(\Gamma)}\otimes
 {\bf H})$ directly by the cap product of $c_{top}(\pi_X^{\ast}
{\bf W}_{canon}^{\circ}|_{Y(\Gamma)}\otimes
 {\bf H})$ with 
 $c_{top}(\pi_X^{\ast}{\bf W}_{canon}/
{\bf W}_{canon}^{\circ}|_{Y(\Gamma)\times T(M)}\otimes {\bf H})$.

\medskip

In this crucial proposition, we observe that the inductively blowup 
procedure modifies the top Chern class of 
$\pi_X^{\ast}{\bf W}_{canon}\otimes {\bf H}$ in a magical way that numerically 
it is possible to identify the top Chern class of 
the modified vector bundle with the cap product of the 
top Chern classes of ${\bf V}_{quot}$ and\footnote{The quotient bundle
 ${\bf V}_{quot}$ of ${\bf W}_{canon}|_{Y(\Gamma)\times T(M)}$ is 
defined in the start of section 3 of [Liu6].} of
$\pi_X^{\ast}{\bf W}_{canon}^{\circ}|_{Y(\Gamma)}
\otimes
 {\bf H}$, the
 canonical algebraic 
obstruction bundle of the new class $C-{\bf M}(E)-\sum_{e_i\cdot (C-{\bf M}(E)E)<0}
e_i$.

The key observation is that the blowing ups along these
 $Z(s_{canon})\cap Y(\Gamma')$, $\Gamma'\in 
\bar{I}_{\Gamma}-\bar{I}_{\Gamma}^{\gg}$, we
\footnote{The index set $\bar{I}_{\Gamma}-\bar{I}_{\Gamma}^{\gg}
\subset \Delta(n)$ consists of the
 admissible graphs $\Gamma'$ such that (i). the strict transforms of
 $Z(s_{canon})\times_{M_{\delta}}Y(\Gamma')$ are blown up ahead of 
 the strict transform of $Z(s_{canon})\times_{M_{\delta}}Y(\Gamma)$ in our initial 
scheme,  (ii). $\Gamma\succ \Gamma'$.
 See defintion 6 on page 47 of [Liu6] for more details.}
 have performed are exactly along the
 sub-loci of $(Z(s_{canon})-Z(s_{canon}^{\circ}))\times_{M_{\delta}}Y(\Gamma)$
 which account for the
 discrepancy $Z(s_{canon})-Z(s_{canon}^{\circ})$ in $Y(\Gamma)$.
 On the other hand, the discrepancy occurs exactly because
 the bundle map $\pi_X^{\ast}{\bf W}_{canon}^{\circ}\otimes {\bf H}
|_{X\times_{M_{\delta}}Y(\Gamma)}
\stackrel{f}{\mapsto} 
\pi_X^{\ast}{\bf W}_{canon}\otimes {\bf H}|_{X\times_{M_{\delta}}Y(\Gamma)}$ 
fails to be injective over some sub-locus of 
$X\times_{M_{\delta}}(Y(\Gamma)-Y_{\Gamma})$. Moreover we actually identify
$\overline{Z(s_{canon})-Z(s_{canon}^{\circ})}\times_{M_{\delta}}Y(\Gamma)$ with 
 the intersection of $s_{canon}^{\circ}$ and the kernel cone of the above 
 bundle map $f$.

 Then in proposition 9 on page 40 of [Liu6], we justify the usage of the blowup
construction to be the right procedure to relate the top Chern classes 
of $\pi_X^{\ast}{\bf W}_{canon}^{\circ}\otimes {\bf H}
|_{X\times_{M_{\delta}}Y(\Gamma)}$ and of $\pi_X^{\ast}
{\bf W}_{canon}\otimes {\bf H}|_{X\times_{M_{\delta}}Y(\Gamma)}$ 

\medskip

\noindent Response to {\bf Problem 6}: The way that we avoid over-counting is
 to demonstrate that our blowup construction is compatible with the 
inclusion-exclusion principle. To illustrate our basic idea, we consider
a toy model regarding the number of elements in a unions of finite sets.

\medskip

Recall a simple fact from the theory of finite sets.
 Let $|A|$ denote the cardinality of the finite set $A$.
 Then $|A\cup B|+|A\cap B|=|A|+|B|$. In other words, we may
 use $|A|+|B|-|A\cap B|$ to calculate $|A\cup B|$.
 This suggests that one may substract the ``over-counting'' term from
 $|A|+|B|$ to get the correct answer.
 On the other hand,
 an alternative way to calculate $|A\cup B|$ is to re-express
 $A\cup B=(A-A\cap B)\cup (A\cap B)\cup (B-A\cap B)$, a disjoint union of
 subsets of $A\cup B$. 
Then $$|A\cup B|=|A-A\cap B|+|A\cap B|+|B-A\cap B|=
 (|A|-|A\cap B|)+|A\cap B|+(|B|-|A\cap B|).$$

 For a finite union of finite sets $A_1, A_2, A_3, \cdots, A_k$, these
two counting schemes lead to two equivalent yet distinct formulae.

 In the first scheme we write 
$|\cup_{1\leq i\leq k}A_i|=
\sum_{I\subset \{1, 2, \cdots, k\}}(-1)^{|I|}|\cap_{i\in I}A_i|$, where
 the index sets run through all the subsets of $\{1, 2, \cdots, k\}$.
 This is the standard form of inclusion/exclusion formula.
 In the second case we define the modified cardinality
 $|\cap_{i\in I}A_i|^{\ast}$ to be 
$|\cap_{i\in I}A_i-\cup_{I\subset J\not=I}\cap_{i\in J}A_i|$, i.e. the
 cardinality of the elements in $\cap_{i\in I}A_i$ which are not
 in any refined intersection $\cap_{i\in J}A_i$ for any $J\supset I$, $J\not=I$.

\begin{lemm}\label{lemm; exclude}
For all $I\subset \{1, 2, \cdots, k\}$, 
let $|\cap_{i\in I}A_i|^{\ast}$ denote the modified cardinality of
 $\cap_{i\in I}A_i$. Then we have the following
 identities relating the modified and the original cardinalities,

$$|\cap_{i\in I}A_i|^{\ast}=|\cap_{i\in I}A_i|-\sum_{J\supset I; J\not=I}
|\cap_{i\in J}A_i|^{\ast},$$

 for all $I\subset \{1, 2, \cdots, k\}$.
\end{lemm}

\noindent Proof: By adjusting the terms in the above identies, 
it suffices to show that for all $I\subset \{1, 2, \cdots, k\}$,

$$|\cap_{i\in I}A_i|^{\ast}+\sum_{J\supset I; J\not=I}
|\cap_{i\in J}A_i|^{\ast}=|\cap_{i\in I}A_i|.$$

 One may write $\cap_{i\in I}A_i$ as the disjoint union
 $$\cap_{i\in I}A_i=
\coprod_{I\subset I'}(\cap_{i\in I'}A_i-\cup_{I'\subset J\not=I'}
(\cap_{i\in J}A_i)),$$

as each element $\in \cap_{i\in I}A_i$ appears in the disjoint union 
on the right hand side
 exactly once.

 By taking cardinalities on both sides and by using the definitions of
 the modified cardinalities, we get the desired equality. 
$\Box$

 In this toy model we ``stratify'' the union $\cup_{i=1}^k A_i$ into
 different set theorectical strata according to the complete collections
 of $A_j$ that an element $\in \cup_{i=1}^k A_i$ belongs to.

 The ``modified'' cardinalities of a finite set counts
the number of elements in a finite intersection of 
$A_j$ which do not lie in the intersection of more $A_j$s. This concept
 is the prototype of the modified family invariants in [Liu1] and [Liu6].

\medskip

 The important characteristics are: 

\medskip

\noindent (a). When there are more and more finite sets involved, the various
intersections of distinct sets lead to a hierachy of intersections and
 a hierachy of ``modified'' cardinalities.

\medskip

\noindent 
(b). To define the modified cardinality numerically, by lemma \ref{lemm; exclude}
 we may start with the
un-modified cardinalities and subtract away all the ``correction terms'' of
 modified cardinalities involving intersections of more $A_j$s.

 It involves an inductive definition. If there are 
 $m$ finite sets $A_i$, $1\leq i\leq k$, all together. Then the modified 
cardinality of $\cap_{i=1}^{i=k}A_i$ is set to 
coincide with the usual (un-modified)
cardinality. By a backward induction decreasing the number of intersections of
 finite
sets $A_j$, at the end one can define the modified cardinalities for all $A_i$.

\medskip

\noindent (c). This alternative approach does not involve the alternating
 sum of terms which are typical to the first formulation of
the inclusion-exclusion principle.

For simplicity, let us illustrate the basic idea how the modified
 algebraic family invariants are defined in [Liu6] 
 by working on the simplified situation $\Gamma_1, \Gamma_2\in 
\Delta(n)-\{\gamma_{\delta}\}$.

The inductive 
blowup construction in section 5 of [Liu6] has been performed in such a way that
when $X$ is blown up along $Z(s_{canon})\times_{M_{\delta}}Y(\Gamma_1)$ and 
along $Z(s_{canon})\times_{M_{\delta}}Y(\Gamma_2)$, a blowing up along
the locus $Z(s_{canon})\times_{M_{\delta}}Y(\Gamma')\supset 
Z(s_{canon})\times_{M_{\delta}}\bigl(Y(\Gamma_1)\cap Y(\Gamma_2)\bigr)$, $\Gamma'
 \in \Delta(n)-\{\gamma_{\delta}\}$ has be done
 in advance. That is why the family invariants we attach to 
 $Y(\Gamma_1)$ or $Y(\Gamma_2)$ get ``modified'' in this process.

 Symbollically let $n_{Y(\Gamma_1)\cup Y(\Gamma_2)}$, $n_{Y(\Gamma_1)}$,
 $n_{Y(\Gamma_2)}$, and $n_{Y(\Gamma')}$ 
denote the intersection numbers

$$\hskip -.3in
\{c_{total}(\pi_X^{\ast}{\bf W}_{canon}\otimes {\bf H}\otimes {\cal O}(-
D_Z))\cap s_{total}(Z, X)\}_{dim_{\bf C}X-rank_{\bf C}{\bf W}_{canon}}
\cap c_1({\bf H})^{rank_{\bf C}({\bf V}_{canon}-{\bf W}_{canon})-1+dim_{\bf C}M_{\delta}}
\cap \{t_L\}$$ built up from the
 localized contributions of top Chern classes attached to (a). 
$Z=Y(\Gamma_1)\cup
 Y(\Gamma_2)$, setting $D_Z=\emptyset$, (b). 
$Z=Y(\Gamma_1)$, or $Z=Y(\Gamma_2)$ and setting $D_Z$ to be the blowup
 exceptional divisor of $Z(s_{canon})\times_{M_{\delta}}Y(\Gamma')$, 
and (c). $Z=Y(\Gamma')$, setting $D_Z=\emptyset$,
 respectively. In this simplified model example the modified
 intersection numbers (called modified algebraic family Seiberg-Witten
 invariants in the long paper [Liu6]) 
  $n_{Y(\Gamma_1)}^{\ast}$,
 $n_{Y(\Gamma_2)}^{\ast}$, and $n_{Y(\Gamma')}^{\ast}$ 
 have to satisfy 

(i). $n_{Y(\Gamma')}^{\ast}=n_{Y(\Gamma')}$.

(ii).$n_{Y(\Gamma_1)}^{\ast}=n_{Y(\Gamma_1)}-
n_{Y(\Gamma')}^{\ast}$.

(iii). $n_{Y(\Gamma_2)}^{\ast}=n_{Y(\Gamma_2)}-
n_{Y(\Gamma')}^{\ast}$.

\medskip

 And finally we have 

$$n_{Y(\Gamma_1)\cup Y(\Gamma_2)}=n_{Y(\Gamma_1)}^{\ast}+
n_{Y(\Gamma_2)}^{\ast}+n_{Y(\Gamma')}^{\ast}.$$

\medskip

 The general cases are dealt with by the same philosophy.
In the following we outline the parallelism between
 the general case and the above toy model.

{\bf I}. The analogue of the inclusion relationship $\subset$ among different 
 finite intersections of the finite sets $A_i$ in family Seiberg-Witten
 theory is coded by
 the partial ordering $\gg$ among the graphs $\Gamma$ 
defined in definition 11 of the long paper [Liu6].

 The partial ordering $\gg$ gives a sufficient condition which guarantees the
subscheme

 \noindent ${\cal M}_{C-{\bf M}(E)E-\sum_{(C-{\bf M}(E)E)\cdot e_i'<0}e_i'}\times_{M_{\delta}}
Y(\Gamma')$ is included in
${\cal M}_{C-{\bf M}(E)E-\sum_{(C-{\bf M}(E)E)\cdot e_i<0}e_i}\times_{M_{\delta}}
Y(\Gamma)=Z(s_{canon}^{\circ})\times_{M_{\delta}}Y(\Gamma)$.

\medskip

{\bf II}. The ``modified algebraic 
family invariants'' are defined following exactly
the same pattern that the ``modified cardinalities'' are defined above.
It involves more technical arguments and the introduction of the classes
 $\tau_{\Gamma}\in K_0(Y(\Gamma)\times T(M))$.

  Starting with some admissible
 graph $\Gamma$ smallest under $\gg$, the modified
 invariant attached to 
 ${\cal M}_{C-{\bf M}(E)E-\sum_{(C-{\bf M}(E)E)\cdot e_i<0}e_i}\times_{M_{\delta}}
Y(\Gamma)$ does not need any correction terms and is set to be
 equal to some mixed family invariant. See definition 12 of [Liu6]
 for more details.

 By an induction argument upon the reversed ordering of 
$\gg$, we define the modified family
 invariants attached to each $Y(\Gamma)$, $\Gamma\in\Delta(n)$.
 This has been done in 
definition 13 and 14 of [Liu6].

\medskip

{\bf III}. In our definitions of modified invariants, 
the alternating sign does not occur. This is because we adopt the
 additive formulation of the inclusion-exclusion principle. It is
 not hard to reformulate our approach into the standard 
inclusion-exclusion formulae involving alternating sums.

\section{The Open Problems Related to the Proof of
Universality Theorem}\label{section; open}

\bigskip

 In this sub-section, we list some open problems related to or inspired from the
 proof of the universality theorem. 

 In [Liu1] we had identified the coefficients of the universal polynomials
 using G${\ddot o}$ttsche' proposal, identifying them with parts of Gromov-Witten
 invariants on ${\bf P}^2$, $K3$, $T^4$ and using 
symplectic technique (including Taubes'
 identification of $SW=Gr$ for smooth pseudo-holomorphic 
curves). On the other hand, the
 construction of residual intersection formula of top Chern classes and 
 the repeatedly scheme-theoretical blowing ups of $X={\bf P}({\bf V}_{canon})$
  in the algebraic category enables us to construct a cycle class
 (the localized top Chern class of some modified algebraic obstruction 
bundle\footnote{We do not write down the modified bundle here. Please
 consult the paper [Liu6] for the detail expression.}) representing the virtual
 fundamental class of moduli 
scheme of curves with $\delta$-node nodal singularities.

 Encouraged by the ${\cal C}^{\infty}$ 
identification of the intersection numbers
 attached to the cycle classes with the explicit enumerations of 
 Gromov-Witten invariants in the concrete examples through pseudo-holomorphic
 curves, we may post the following open question:

\noindent {\bf Open Question}: 
Identify the generating functions of the universal polynomials algebraically.

The universality theorem asserts the existence
 of these universal polynomials for all $\delta\in {\bf N}$ which code the
 enumerative information on the ``number of nodal curves''. By an ingenious
 argument of G${\ddot o}$ttsche he has shown [Got] that the generating function of these
 universal polynomials takes a factorizable form 

 $${\cal F}(q)=
A_1(q)^{c_1(K_M)^2} A_2(q)^{c_2(M)} A_3(q)^{c_1(L)^2} A_4(q)^{c_1(L)\cdot
 c_1(K_M)}.$$

Substituting $q$ by\footnote{$D=q{d\over dq}$ and $G_2(q)$ 
is the quasi-modular form ${-1\over 24}+\sum_{k\geq 1}\sigma_1(k)q^k$.}
 $DG_2(q)$, then ${\cal F}(DG_2(q))$ can be identified with
 $${(DG_2(q)/q)^{\chi(L)}B_1(q)^{c_1(K_M)^2}B_2(q)^{c_1(L)\cdot c_1(K_M)}\over
 (\Delta(q)D^2G_2(q)/q^2)^{\chi({\cal O}_M)\over 2}},$$ 
where $B_1$ and $B_2$ are the two power series derived by G{$\ddot o$}ttsche
[Got] starting with $$B_1(q)=1-q-5q^2+30q^3-345q^4+2961 q^5\ldots,$$
  and $$B_2(q)=1+5q+2q^2+35q^3-140q^4+986q^5+\ldots.$$

 On the other hand, the original ${\cal C}^{\infty}$ approach of the
 problem suggests a symplectic generalization of the universality theorem.

\noindent {\bf Open Question}: Formulate and 
Generalize the universality theorem to
 symplectic four-manifolds $M$ with/without the $b^+_2=1$ condition.

\medskip

\noindent {\bf Open Sub-Question}: For symplectic four-manifolds with
 $b^+_2>1$ (which correspond to algebraic surfaces with $p_g>0$), 
 find the right analogue for the algebraic family Seiberg-Witten invariants
 which correspond non-trivial family Seiberg-Witten invariants.

\medskip

\noindent {\bf Open Sub-Question}:
 Find the suitable definition of ``number of nodal curves'' or
 more generally the ``number of singular curves'' in the symplectic
 category. Our study based 
on algebraic geometric means suggests that one may need to work in the framework
of the ``virtual numbers'' instead of the discrete number count of curves 
in a given class $\in H^2(M, {\bf Z})$.

 The residual obstruction vector bundle 
construction in the proof of the universality theorem not only defines 
 the integer valued ``number of nodal curves'', it also constructs the 
virtual fundamental class of $\delta-$node singular curves in a 
``$5\delta-1$''-very ample linear system. On the other hand, 
 the technique of residual intersection theory upon type $II$ exceptional
 curves [Liu7] allows us to drop the $5\delta-1$-very ample condition and 
define the virtual fundamental class of the moduli space of
 nodal curves for general linear systems on general algebraic surfaces.

  In our foundation the nodal curve invariant counts embedded nodal curves. 
On the other hand, the Gromov-Witten theory counts maps, including the various
 types of 
multiple covering maps induced from different fractions of the class. 

Unlike Gromov-Witten theory which are ${\bf Q}$ valued, 
our theory does not make use of the moduli stacks of curves.
 Thus our invariants are naturally ${\bf Z}$ valued. This fits to the
 philosophy of Gopakumar-Vafa conjecture [GV], [BP], [HST] in constructing
 ${\bf Z}$-valued curve counting invariants. 
 
Thus, it is very desirable to find out:

\medskip

\noindent {\bf Open Question}: Multiple Covering Formula--Find out the
weighted contribution of the multiple covering maps and relate our
``nodal curve'' invariant with the usual Gromov-Witten invariant when
 $p_g=0$. When $p_g=1$ and $M=K3$ or $T^4$, we may consider the
 ${\bf S}^2$ hyperkahler families of $K3$ and the ${\bf S}^2$ family
 Gromov-Witten invariant and formulate a similar problem. 

\medskip

Another interesting question worthy to pursue is,

\medskip

\noindent {\bf Open Question}: Go beyond the range of nodal curve singularities
 and study the structure of ``number of singular curves'' with non-nodal 
singularities. 
 
 In particular, the algebraic proof of the universality theorem along with the
 finiteness result of the curves $\in |L|$
for an effective bound on the very-ampleness of $L$ also
 work for curves with ordinary singularities with multiplicity $m>2$.

\medskip

\noindent {\bf Open Sub-Question}: Generalize the Caporaso-Harris [CH]
 recursive formula of counting of nodal curves on ${\bf P}^2$ 
to curves with ordinary singularities with multiplicity $m>2$.

\medskip

 Caporaso-Harris original argument could be intrepreted naturally in
 terms of Gromov-Witten invariants of ${\bf P}^2$ (see e.g. section 7 of 
[Va]), 
the fact that their
 original argument does not rely on the usage of Gromov-Witten invariants
\footnote{for non-nodal curves, it is unclear at this moment that there are  
 ``Gromov-Witten type invariants'' corresponding to them.}
 may look encouraging to us.

\medskip

\noindent {\bf Open Problem}: The local proof of the Blowup formula.

 If one blows up the algebraic surface $M$ at a single point and get $\tilde{M}$, 
then the pull-back of a $5\delta-1$-very ample line 
 bundle $L$ on $M$ fails to be very ample on the blown up surface  $\tilde{M}$ 
exactly
 at the exceptional locus of $\tilde{M}\mapsto M$. 
On the other hand, a simple calculation
 on G${\ddot o}$ttsche's formula indicates that the formula of nodal curves changes
 by the following universal formula (again substituting $q$ by $DG_2(q)$). 

$${\cal F}_{\tilde M}(DG_2)={\cal F}_{M}(DG_2)\cdot
 ({B_2(q)\over B_1(q)})\cdot ({DG_2\over q})^{-1}.$$

The power series $B_1(q)$, $B_2(q)$ are determined from J. Harris and
Caperaso's calulation of Severi degrees 
[CH] by G${\ddot o}$ttsche [Got] modulo Yau-Zaslow formula, etc.

The blowup formula implies that a purely ``local way''
 to identify the coefficients
 of $B_1(q), B_2(q)$ is possible. 
In fact, such an identification will relate intersection
 numbers of local nature to the global enumerative geometry 
datum gethered from Caporaso-Harris
 calculation [CH] and it is very interesting to understand their
 relationship.

\medskip

\noindent {\bf Open Sub-Problem}: Understand the intrinsic geometric meaning
 of the power series $B_1(q)$ and $B_2(q)$.

 By comparing the Yau-Zaslow formula ${1\over \prod_i(1-q^i)}^{c_2(M)}$, $M=K3$
 with the Riemann-Roch formula, it is clear that ${1\over \prod_i(1-q^i)}$ is
 the natural prolongation of the coefficient ${1\over 12}$ in front of $c_2(M)$ 
in the classical Noether formula.
It is desirable to understand the relationship of $B_1(q)$ and $B_2(q)$ with
 the other coefficients of the surface Riemann-Roch formula.

 Our formulation of the universality theorem
as the natural prolongation of Riemann-Roch theorem suggests
 that the generating power series of ``numbers of singular curves'' 
 with different prescribed topological types of curve singularities
 should be related to each other and prolong the classical Riemann-Roch
 formula in a natural way. So we propose to
 study the topological structure of curve singularities
 in order to find recursive formulae between different singularities.

 It is known that the topology of an isolated algebraic 
curve singularity is completely
 coded by the topology of the link space $\subset {\bf S}^3$ of the
 singularity. Then the topology
 of the link is coded by the Alexander polynomial of the link.

 Thus, one may ask the following question,

\noindent {\bf Open Question}: Find out the dependence of the
 universality formulae on the alexander polynomials of the link spaces of
the singularities.

Evidence of the nodal curve case shows that a part of the generating function
 has been modular on Calabi-Yau algebraic surfaces (i.e. $K3$ or $T^4$).

 We expect to get liftings from Axelander polynomials of links
 to modular objects in this generalization.

  We also mention that by generalizing the tools we have used to
 higher dimensions, the universality theorem can be generalized to
 higher dimensions. Instead of counting curves with prescribed singularities,
  we enumerate divisors with prescribed singularities in very
  ample linear systems. Then the program involves prolonging the whole
 $Todd_M\cdot ch(L)$ expression into multiplicative power series of powers of
 the Chern numbers $c_1^{a_1}(M)c_2^{a_2}(M)\cdots c_n^{a_n}(M)c_1^k(L)$, 
 with $\sum ia_i+k=n$.

\bigskip

 Finally, we have noticed earlier (see page \pageref{DT}) 
that the simple invariants that we have attached to linear systems are the 
degenerated cases (reduced to $dim_{\bf C}M=2$) 
of Donaldson-Thomas invariants of Hilbert schemes of curves.
 Our machinery has provided a tower of enumerative invariants of singular
 curves from the leading invariant of linear systems.

\medskip

\noindent {\bf Open Question}: Generalize the concept of Nodal Curve Invariants
 to the case of Calabi-Yau three-folds. The algebraic construction of nodal
 curve invariants is supposed to
 probe the singular curves inside the Hilbert schemes of
 curves such that Donaldson-Thomas invariants are the first order leading 
invariants. It should be closed related to the hypothetical integral valued 
Gopakumar-Vafa invariants proposed in [GV].

{}
\end{document}